\documentclass[a4paper,12pt]{article}

\usepackage{emptypage} 
\usepackage{a4wide}

\usepackage[pagestyles]{titlesec}
\usepackage{amsmath}
\usepackage{accents}

\newcommand{\ubar}[1]{\underaccent{\bar}{#1}}

\newpagestyle{main}{
\headrule
\sethead[\thepage][\chaptertitle][] 
{}{\sectiontitle}{\thepage} 
}

\usepackage{latexsym}
\usepackage{amssymb}
\usepackage{amsthm}

\usepackage{graphicx}

\usepackage{amstext}
\usepackage{bbm}

\usepackage{xcolor}
\usepackage{amsmath}
\usepackage{hyperref}
 
\usepackage[capitalize,nameinlink]{cleveref}
\usepackage{autonum}

\usepackage[]{currvita}

\usepackage[math]{cellspace}
\makeatletter
\def\underbracex#1#2{\mathop{\vtop{\m@th\ialign{##\crcr
   $\hfil\displaystyle{#2}\hfil$\crcr
   \noalign{\kern3\p@\nointerlineskip}%
   #1\crcr\noalign{\kern3\p@}}}}\limits}

\def\upbracefilla{$\m@th \setbox\z@\hbox{$\braceld$}%
  \bracelu\leaders\vrule \@height\ht\z@ \@depth\z@\hfill 
\kern\p@\vrule \@width\p@\kern\p@\vrule \@width\p@\kern\p@\vrule \@width\p@
$}

\def\upbracefillb{$\m@th \setbox\z@\hbox{$\braceld$}%
\vrule \@width\p@\kern\p@\vrule \@width\p@\kern\p@\vrule \@width\p@\kern\p@
 \leaders\vrule \@height\ht\z@ \@depth\z@\hfill\bracerd
  \braceld\leaders\vrule \@height\ht\z@ \@depth\z@\hfill
\kern\p@\vrule \@width\p@\kern\p@\vrule \@width\p@\kern\p@\vrule \@width\p@
$}

\def\upbracefillc{$\m@th \setbox\z@\hbox{$\braceld$}%
\vrule \@width\p@\kern\p@\vrule \@width\p@\kern\p@\vrule \@width\p@\kern\p@
\leaders\vrule \@height\ht\z@ \@depth\z@\hfill
\kern\p@\vrule \@width\p@\kern\p@\vrule \@width\p@\kern\p@\vrule \@width\p@
$}

\def\upbracefill{$\m@th \setbox\z@\hbox{$\braceld$}%
\vrule \@width\p@\kern\p@\vrule \@width\p@\kern\p@\vrule \@width\p@\kern\p@
 \leaders\vrule \@height\ht\z@ \@depth\z@\hfill\braceru$}

\def\upbracefillL{$\m@th \setbox\z@\hbox{$\braceld$}%
  \bracelu\leaders\vrule \@height\ht\z@ \@depth\z@\hfill 
\kern\p@\vrule \@width\p@\kern\p@\vrule \@width\p@\kern\p@\vrule \@width\p@
$}

\def\upbracefillR{$\m@th \setbox\z@\hbox{$\braceld$}%
\vrule \@width\p@\kern\p@\vrule \@width\p@\kern\p@\vrule \@width\p@\kern\p@
 \leaders\vrule \@height\ht\z@ \@depth\z@\hfill\bracerd
  \braceld\leaders\vrule \@height\ht\z@ \@depth\z@\hfill
 \leaders\vrule \@height\ht\z@ \@depth\z@\hfill\braceru$}

\def\BRR{\underbracex\upbracefillLR}

\def\upbracefillLR{$\m@th \setbox\z@\hbox{$\braceld$}%
  \bracelu\leaders\vrule \@height\ht\z@ \@depth\z@\hfill 
\leaders\vrule \@height\ht\z@ \@depth\z@\hfill\bracerd
  \braceld\leaders\vrule \@height\ht\z@ \@depth\z@\hfill
 \leaders\vrule \@height\ht\z@ \@depth\z@\hfill\braceru
$}
\makeatother
\newcommand{\retainlabel}[1]{\label{#1}\sbox0{\ref{#1}}}

\def\R{\mathbb R}

\def\N{\mathbb N}

\newcommand{\p}{\psi}
\renewcommand{\t}{\theta}

\newcommand{\ti}{\tilde }

\newcommand{\hatp}{{\hat\p}}
\newcommand{\hatt}{{\hat \t}}

\newcommand{\hxi}{\hat \xi}
\newcommand{\hu}{\hat \u}

\newcommand{\ten}{{\t_n^{\eps}}}
\newcommand{\pen}{\p_n^{\eps}}

\newcommand{\luen}{{\lambda_{n}^{\eps}}}

\def\u{u}
\def\g{\gamma}

\def\a{\alpha}

\renewcommand{\v}{v}
\newcommand{\w}{w}

\newcommand{\pa}{\partial_}

\newcommand{\ds}{\fr d {ds}}

\newcommand{\dxi}{\pa \xi}

\newcommand{\du}{\pa \u}

\newcommand{\dx}{\pa x}

\newcommand{\deps}{\pa \eps}
\def\xx{_{xx}}

\renewcommand{\d}{\dot}

\def \F{F(\eps,x)}
\def \tiF{\ti F(\eps)}

\def \Go{{\cal G}_1^{\eps}}
\def \Gt{{\cal G}_2^{\eps}}

\def \Gn{{\cal G}_n^{\eps}}

\def \tiG{\tilde{\cal G}^{\eps}}

\def \tiGo{\tilde{\cal G}_1^{\eps}}
\def \tiGt{\tilde{\cal G}_2^{\eps}}

\def \tiGn{\tilde{\cal G}_n^{\eps}}

\newcommand{\BR}[2]{\BRR{{#1}}_\text{{#2}}}

\newcommand{\OT}[2]{\stackrel{{\smash{\scriptscriptstyle{#2}}}}{{#1}}}

\newcommand{\Ltwortwoaxix}[2]{\l\langle {#1},{#2} \r\rangle_{L^{2,a}_{\xi,x}(\R^2)\oplus L^{2,a}_{\xi,x}(\R^2)}}

\newcommand{\nw}[1]{\l|{#1}\r|_{L^2(\R)}}

\newcommand{\nltwo}[1]{\l|{#1}\r|_{L^2(\R)}}

\newcommand{\nhone}[1]{\l|{#1}\r|_{H^1(\R)}}

\newcommand{\fr}{\frac}

\def\be#1\ee{\begin{align}#1\end{align}}
\def\ben#1\een{\begin{align+}#1\end{align+}}
\def\ba#1\ea{\begin{aligned}[t]#1\end{aligned}}
\def\bs#1\es{\begin{split}#1\end{split}}

\renewcommand{\l}{\left}
\renewcommand{\r}{\right}
\def\?{????}

\newcommand{\bma}{\begin{pmatrix}}
\newcommand{\ema}{\end{pmatrix}}

\newcommand{\bmat}{\begin{bmatrix}}
\newcommand{\emat}{\end{bmatrix}}

\newcommand{\bca}{\begin{cases}}
\newcommand{\eca}{\end{cases}}

\newcommand{\re}{\eqref}
\newcommand{\nn}{\nonumber}

\newcommand{\la}{\label}
\newcommand{\beq}{\begin{equation}}
\newcommand{\eeq}{\end{equation}}
\newcommand{\eps}{\varepsilon}

\renewcommand{\qed}{\protect~\protect\hfill $\Box$}

\begin{document}


\newtheoremstyle{plainNEW}
  {}
  {}
  {\itshape}
  {}
  {\boldmath\bfseries}
  {.}
  { }
  {\thmname{#1}\thmnumber{ #2}\thmnote{ (#3)}}%


\theoremstyle{plainNEW}

\newtheorem{theorem}{Theorem}[section]

\newtheorem{definition}[theorem]{Definition}
\newtheorem{deflem}[theorem]{Definition and Lemma}
\newtheorem{lemma}[theorem]{Lemma}
\newtheorem{corollary}[theorem]{Corollary}
\newcommand{\bde}{\begin{definition}}
\newcommand{\ede}{\end{definition}}
\newcommand{\ble}{\begin{lemma}}
\newcommand{\ele}{\end{lemma}}
\newcommand{\bre}{\begin{remark}}
\newcommand{\ere}{\end{remark}}
\newcommand{\bco}{\begin{corollary}}
\newcommand{\eco}{\end{corollary}}
\newcommand{\bpro}{\begin{proposition}}
\newcommand{\epro}{\end{proposition}}

\newtheorem{assumption}{Assumptions}[theorem]
\newcommand{\bas}{\begin{assumption}}
\newcommand{\eas}{\end{assumption}}

\theoremstyle{definition}

\newtheorem{example}[theorem]{Example}
\newtheorem{remark}[theorem]{Remark}
\newtheorem{remarks}[theorem]{Remarks}
\newtheorem*{proofNEW}{Proof}
\newcommand{\bpr}{\begin{proofNEW}}
\newcommand{\epr}{\qed
\bigskip
\end{proofNEW}}

\theoremstyle{plainNEW}

\newcommand{\bth}{\begin{theorem}}
\renewcommand{\eth}{\end{theorem}}
\newtheorem{proposition}[theorem]{Proposition}

\pagenumbering{gobble}
\newpage
\title{Invariant Virtual Solitary Manifold\\ 
of the Perturbed Sine-Gordon Equation }

\date{}

\author{{\sc Timur Mashkin}\\[2ex]
         Mathematisches Institut, Universit\"at K\"oln, \\
         Weyertal 86-90, D\,-\,50931 K\"oln, Germany \\
         e-mail: tmashkin@math.uni-koeln.de}

\maketitle

\begin{abstract}\noindent
We study the perturbed sine-Gordon equation 
$\theta_{tt}-\theta_{xx}+\sin \theta=  F(\eps,x)$,
where
we assume that 
the perturbation $F$ is analytic in $\eps$ and that its derivatives with respect to $\eps$ satisfy certain bounds at $\eps=0$. 
We construct implicitly an, adjusted to the perturbation $F$, virtual solitary manifold, which is invariant in the following sense: 
The initial value problem for the perturbed sine-Gordon equation with an appropriate initial state on the constructed manifold has a unique solution, which follows a trajectory on the
virtual solitary 
manifold. The trajectory is precisely described by two parameters, which satisfy a specific system of ODEs. 

The approach is based on \cite{Mashkin}, where we constructed by an iteration scheme a virtual solitary manifold for the perturbed sine-Gordon equation.
In \cite{Mashkin} we proved a stability result for the perturbed sine-Gordon equation with initial data close to the virtual solitary manifold.
The employed iteration scheme produces a sequence of virtual solitary manifolds such that the accuracy of the corresponding stability statements  
increases after each iteration step, as long as the perturbation $F$ is sufficiently often differentiable.
The invariant virtual solitary manifold constructed in this work is generated as a limit of 
the virtual solitary manifolds produced by the iteration scheme. 

The method and the kind of result presented in this paper is to our knowledge a novelty in the field of stability of solitons.
\end{abstract}

\pagenumbering{arabic}

\section{Introduction}

\noindent 
We consider the perturbed sine-Gordon equation 
\beq\la{SGE}
\t_{tt}-\t_{xx}+\sin\t=F(\eps,x),~~~~t,x\in\R,~~~~\eps\ll 1,
\eeq
which can be written as a system in first order formulation:
\be\la{SGE1 first order introduction}
\partial_t\bma
\t\\
\p
\ema=\l(\begin{matrix}
\p\\
\t_{xx}-\sin\t+F(\eps,x)\\
\end{matrix}\r).
\ee
The unperturbed sine-Gordon equation (i.e., $F(\eps,x)=0$)
admits soliton solutions 
$$
\bma 
\t_0(\xi(t),\u(t),x)\\
\p_0(\xi(t),\u(t),x)
\ema,
~\text{where}
$$
\be\la{ODEintro}
\dot\xi=\u\,, ~~ \dot \u=0\,,~~~~(\xi(0),\u(0))=(a,v)\in\R\times(-1,1).
\ee
Here the functions $(\t_0,\p_0)$ are defined by
\be\la{solitonsolution}
{}&\bma
\t_0(\xi,\u,x)\\
\p_0(\xi,\u,x)
\ema:=\bma
\t_K(\g(u)(x-\xi))\\
-\u\g(u)\t_K'(\g(u)(x-\xi))\\
\ema\,,~\u\in(-1,1),~~\xi,x\in\R,
\ee
where 
$$
\g(u)=\frac 1 {\sqrt{1-u^2}},~~~~\t_K(x) =4\arctan(e^x),
$$
and $\t_K$ satisfies $\t_K''(x)=\sin\t_K(x)$ with boundary conditions $\t_K(x) \to \bma 2\pi\\0  \ema$ as $x\to \pm \infty$.
The states $\l(\begin{matrix}
\t_0(a,v,\cdot)\\
\p_0(a,v,\cdot)\\
\end{matrix}\r)$ form the two dimensional classical solitary manifold  
$$
{\cal S}_0:=\l\{ \l(\begin{matrix}
\t_0(a,v,\cdot)\\
\p_0(a,v,\cdot)\\
\end{matrix}\r)~:~v\in(-1,1),~a\in\R
\r\}.
$$

Let us mention some previous works before we state the main result. 
Orbital stability of soliton solutions under perturbations of the initial data has been proven for the unperturbed sine-Gordon equation 
(see \cite{MR678151}, \cite[Section 4]{Stuart3}).
D. M. Stuart \cite{Stuart2} 
considered the perturbed sine-Gordon equation
\be
{}&\t_{tt}-\t_{xx}+\sin\t+\eps g=0,\nn
\ee
for specific perturbations of the form 
$
g=g(\eps t, \eps x, \t)
$ and initial data $\eps$-close to
a kink.
He proved the existence of solutions, which approximate kinks with slowly evolving in time  centre and velocity, up to time $1/\eps$ and up to errors of order $\eps$.
Kinks are solutions of the unperturbed  equation  \re{SGE}, 
given by 
$
\t(t,x)=\t_0(\xi(t),\u(t),x)
$,
where the centre $\xi$ and the velocity $\u$ satisfy ODEs \re{ODEintro}.
The proof is based on an orthogonal
decomposition of the solution into an oscillatory part and a one-dimensional
"zero-mode" term.

In \cite[Part I]{MashkinDissertation} we studied equation \re{SGE1 first order introduction} for different types of perturbations. 
For instance, we proved for $F(\eps,x)= \eps f(\eps x)$ 
that the Cauchy problem for initial data $\eps^\fr{1}2$-close to the classical solitary manifold  
$
{\cal S}_0  
$
has a unique solution, which follows up to time $1/\eps^{\fr 1 4}$ and errors of order $\eps^\fr{1}2 $ a trajectory
on 
$
{\cal S}_0  
$, where the trajectory on $
{\cal S}_0  
$ is described precisely by ODEs for uniform linear motion.  
One should take into account that our perturbation $F(\eps,x)= \eps f(\eps x)$ is not comparable to the perturbations in \cite{Stuart2} due to some specific assumptions made on $g$.

For perturbations of type $F(\eps,x)= \eps^2 f(\eps x)$ with 
$f\in H^3(\R)$, we obtined richer dynamics on the solitary manifold in \cite{MashkinElectricField}. We proved that the Cauchy problem for initial data $\eps^\fr{11}8$-close to the classical solitary manifold  
$
{\cal S}_0  
$
has a unique solution, which follows up to time $1/\eps$ and errors of order $\eps^\fr{3}4$ a trajectory
on 
$
{\cal S}_0  
$. 
The trajectory
on 
$
{\cal S}_0  
$
is described precisely by 
ODEs, 
which contain the perturbation $f$.
The ODEs are obtained by considering restricted Hamilton equations and describe a 
fixed nontrivial perturbation of the uniform linear motion as $\eps \to 0$ if $f(0)\not=0$. 
The evolution of the dynamics on the solitary manifold in \cite[Part I]{MashkinDissertation}/ \cite{MashkinElectricField} is described more accurate than the evolution of the approximated kink in \cite{Stuart2} in the following sense: In \cite[Part I]{MashkinDissertation}/ \cite{MashkinElectricField} the parameters of the manifold satisfy exactly specific ODEs, whereas in \cite{Stuart2} the evolution of the kink parameters are determined just up to errors of order $\eps$.

The proofs of \cite[Part I]{MashkinDissertation}, \cite{MashkinElectricField}, and \cite[Section 4]{Stuart3} are based on a nowadays conventional method for verification of stability of solitons (for different equations), namely the decomposition of the dynamics into a part on the classical solitary manifold and a transversal part along with the application of Lyapunov-type arguments. 
This approach emerges, for instance, also in \cite{MR2094474,MR2232367,MR2342704,HoZwSolitonint,MR2855072}.

In \cite{Mashkin} we extended this method 
by utilizing a virtual solitary manifold. 
There we studied the sine-Gordon equation with perturbations
$
\eps \mapsto F(\eps,\cdot)$
of class
$
C^{n}
$ 
(mapping into a specific weighted Sobolev space on $\R$),
whose first $k$ derivatives vanish at 0, i.e., 
$\deps^l F(0,\cdot)=0~~ \text{for} ~~0\le l\le k$, where $k+1 \le n$ and $n\ge 1$.
We constructed in \cite{Mashkin} by an iteration scheme composed of $n$ steps
a virtual solitary manifold, which is adjusted
to the perturbation $F $. The iteration process
can be thought of as a stepwise distortion of the classical solitary manifold ${\cal S}_0$.
Each step in the iteration scheme corresponds to solving implicitly a specific PDE. 
The implicit solution  $\eps\mapsto (\t_n^\eps(\xi,u,x),
\p_n^\eps(\xi,u,x),\lambda_{n}^\eps\l(\xi, \u\r))$ obtained in
the last iteration step defines
the 
virtual solitary manifold
\be\la{virtMF}
{\cal S}_n^\eps:=\l\{ \bma
\t_n^\eps(a,v,\cdot)\\
\p_n^\eps(a,v,\cdot)
\ema~:~v\in(-u_*,u_*),~a\in\R
\r\}, ~~~~u_*\in (0,1],
\ee
and is used to formulate
the result of \cite{Mashkin}, which is as following:
%
For $\xi_s\in\R$, $\eps\ll 1$,
the Cauchy problem
\be\la{Cauchy_intro}
\partial_t \bma
\t \\
\p 
\ema{}&=\l(\begin{matrix}
\p \\
\dx^2\t -\sin\t +\F\\
\end{matrix}\r),~
\bma
\t(0,x)\\
\p(0,x)
\ema{}&=
\bma
\t^\eps_n(\xi_s,\u_s,x)\\
\p^\eps_n(\xi_s,\u_s,x)
\ema+
\bma
\v(0,x)\\
\w(0,x) 
\ema,
\ee
with appropriate initial data that is $\eps^n$-close to ${\cal S}_n^\eps$, i.e., 
$
\nhone{v(0,\cdot)}^2+\nw{w(0,\cdot)}^2\le \eps^{2n},
$ 
with initial velocity that satisfies the smallness assumption
$
|u_s|\le \tilde C\eps^{\fr{k+1}2} 
$,
has a unique solution $(\t,\p)$,  
which may be written up to time
$ 1/ (\tilde C {\eps} ^{\fr{k+1}2})$ in the form
\be
\bma
\t(t,x)\\ 
\p(t,x)
\ema=\bma 
\t_n^\eps(\bar\xi(t),\bar\u(t),x)\\ 
\p_n^\eps(\bar\xi(t),\bar\u(t),x)
\ema+ 
\bma
\v(t,x)\\
\w(t,x)
\ema.\nn
\ee
The solution remains $\eps^n$-close to ${\cal S}_n^\eps$, i.e.,
$
\nhone{v(t,\cdot)}^2+\nltwo{w(t,\cdot)}^2\le \tilde C \eps^{2n} ,
$
and the dynamics on ${\cal S}_n^\eps$ is described precisely by the parameters $(\bar\xi(t),\bar\u(t))$, which satisfy exactly the ODEs
\be\la{ODE introduction}
 \d{\bar\xi}( t) =  \bar\u(t)  \,,~~~~
 \d{\bar\u}( t) = \lambda_{n}^\eps\l(\bar\xi(t), \bar\u(t)\r), 
\ee
with initial data 
$
\bar\xi(0)=\xi_s,~\bar\u(0)=\u_s
$. 
The parameters $\bar\xi,\bar\u$ describe
a fixed nontrivial perturbation of the uniform linear motion
as $\eps \to 0$ if the perturbation $F$ satisfies a specific condition. 
The higher the differentiability class $C^n$
of $F$ the higher is the accuracy of the stability statement and the more first derivatives of $F$ vanish at 0 
the larger is the time scale of the result.

The sine-Gordon equation arises in various physical applications presented for instance in 
\cite{0953-8984-7-2-013,RevModPhys.61.763,FrenKont,0022-3719-11-1-007}.
In \cite{Skyrme237} T. H. R. Skyrme proposed 
the equation to model elementary particles  and in \cite{doi:10.1143/JPSJ.46.1594}
dynamics of solitons under constant electric field were examined numerically.
We focus in the present work, as also in \cite{Mashkin},
on the interaction  
of virtual solitons 
with a time independent electric field $F(\eps,x)$, which is a physically relevant problem.
\paragraph{Main Result and Consequences}
The iteration scheme introduced in \cite{Mashkin} provides a sequence of implicitly given functions. 
In the present paper, we show that under some additional assumptions the provided sequence, 
denoted by $(\t^\eps_n,\p^\eps_n,\lambda^\eps_n)$,
converges to a limit, which we denote by $(\t^\eps_\infty,\p^\eps_\infty,\lambda^\eps_\infty)$. 
Our main result states that the 
virtual solitary manifold defined analogously to \re{virtMF} by the 
functions
$(\t^\eps_\infty,\p^\eps_\infty,\lambda^\eps_\infty)$
is invariant. 
In greater detail, the main result is as follows. 
Assume that the perturbation 
$
\eps \mapsto F(\eps,\cdot)
$ 
is analytic (mapping into a specific weighted Sobolev space on $\R$), where the derivatives with respect to $\eps$ of $F$ satisfy specific bounds at $\eps=0$ (stated below in \re{assumption bounds derivatives F}) and $  F(0,\cdot)=0$, $\deps  F(0,\cdot)=0$. 
Let $\xi_s\in \R $ and consider the Cauchy problem
\be \la{Cauchy_intro InfMF}
\partial_t \bma
\t \\
\p 
\ema
=\l(\begin{matrix}
\p \\
\dx^2\t -\sin\t +\F\\
\end{matrix}\r) ,~~
%
\bma
\t(0,x)\\
\p(0,x)
\ema=\bma
\t^\eps_\infty(\xi_s,\u_s,x)\\
\p^\eps_\infty(\xi_s,\u_s,x)
\ema ,~~
\eps\ll 1,
\ee 
where the initial velocity satisfies the assumption $|u_s|< \u_*$ for a specific $\u_*$.
Then the Cauchy problem
\re{Cauchy_intro InfMF}
has a unique solution, which may be written in the form
\be\la{intro form ivariant mf} 
\bma
\t(t,x)\\ 
\p(t,x)
\ema=\bma 
\t_\infty^\eps(\bar\xi(t),\bar\u(t),x)\\ 
\p_\infty^\eps(\bar\xi(t),\bar\u(t),x)
\ema,
\ee
where the parameters $(\bar\xi(t),\bar\u(t))$ satisfy the ODEs
\be \la{intro ODEs invariant mf}
 {}&\d{\bar\xi}( t) =  \bar\u(t)  ,~~
\d{\bar\u}( t) = \lambda_{\infty}^\eps\l(\bar\xi(t), \bar\u(t)\r),
\ee
with initial data 
$
\bar\xi(0)=\xi_s,~\bar\u(0)=\u_s 
$.
The solution exists and has this form as long as the parameters stay in an appropriate pareameter area, i.e., as long as $|\bar\xi(t)| \le \Xi, ~|\bar \u(t)|<\u_*$, where $\Xi$ depends on the initial centre $\xi_s$.
In particular, if $|u_s|\le \tilde C \eps^{}$ for a specific $\tilde C$,
then the unique solution exists and can be expressed in the presented form on the time scale
\be\la{intro time scale}
0\le t \le \fr {1} {\tilde C \eps }.
\ee
If additionally the perturbation $F$ satisfies condition \re{intro condition on nontrivial dynamic} mentioned below, then  
the 
parameters $\bar\xi,\bar\u$ describe,
on the nontrivial time scale \re{intro time scale},
a fixed nontrivial perturbation of the uniform linear motion 
as $\eps \to 0$.

The result states that the solution remains on the virtual solitary manifold 
defined by $(\t^\eps_\infty,\p^\eps_\infty)$ 
and it yields a precise description of the solution $(\t,\p)$ to the Cauchy problem \re{Cauchy_intro InfMF}, since the dynamics on the 
manifold is exactly characterized by the ODEs \re{intro ODEs invariant mf}. 
The maximal interval of existence (time interval) of the solution depends on the perturbation $F$ and on the initial data, which determine the ODEs \re{intro ODEs invariant mf}, whereas the ODEs determine for how long the parameters $(\bar\xi(t),\bar\u(t))$ stay in the corresponding parameter area.
A precise statement is found in \cref{se: Main Results}. 
 
The existence of the invariant virtual solitary manifold has a tremendous theoretical value. Furthermore, the invariant manifold allows us to describe the solution of \re{SGE1 first order introduction} with appropriate initial data by far more accurate than it was done in \cite{Mashkin}.
Our main result can be considered as an extension of the work of \cite{Mashkin}, where we corrected the classical solitary manifold of the sine-Gordon equation  
arbitrarily many times (finite number) and 
improved the accuracy of the stability statement in each correction step. 
In this paper the invariant virtual solitary manifold is generated by a limit process - that is, in infinitely many correction steps - in such a way that the manifold is adjusted to the perturbation term $F$.
There exists a 
community, 
which advocates
the following conjecture
for specific PDEs with soliton solutions: For appropriate classes of solutions to the corresponding PDE
there exists a manifold, which acts as an attractor. One expects that for appropriate initial data, not necessarily close to the manifold, the solution is going to come close to the manifold for advancing times. In case of the sine-Gordon equation the virtual solitary manifold generated in this paper is a serious candidate for such an attractive manifold, which makes our result even more interesting for further investigations. 

Our approach and the fact of existence of an invariant manifold for an integrable equation with an external perturbation (invariant in the sense of our main result), 
is to our knowledge a novelty in the field of stability of solitons.
However, singular corrections of the classical solitary manifold have been carried out in other works in different forms such as in \cite{HolmerLin} and in \cite{HoZwSolitonint} for the NLS equation, which corresponds to the first iteration in the scheme from \cite{Mashkin}. 
The idea of modifying the classical solitary manifold of the sine-Gordon equation by utilizing implicitly defined functions appears in \cite[Section 3]{Stuart3}, where the purpose was to rewrite the Hamiltonian in a neighbourhood of the manifold of virtual solitons. 
Neither the virtual solitary manifold \re{virtMF} nor the iteration scheme introduced in \cite{Mashkin} were considered in \cite{Stuart3}.

Several long (but finite)-time results for different 
equations with external potentials can be found, for example, in \cite{MR2094474,MR2232367,MR2342704,MR2855072}.
Further results on orbital stability and long time soliton asymptotics are presented in 
\cite{MR820338,MR0428914,MR0386438,MR2920823,MR1071238,MR1221351,ImaykinKomechVainberg,MR3630087,MR3461359}.
\paragraph{Our Techniques}
We generate the invariant virtual solitary manifold
by utilizing 
the iteration scheme from \cite{Mashkin}, 
whereby we modify the scheme in certain points.
In the present paper, the scheme is
implemented for an analytic function  
$\eps\mapsto \tilde F(\eps)$
mapping into a specific 
Sobolev space on $\R^2$
such that $\tilde F(\eps)$ depends on $(\xi,x)$ (for the sake of clarity, we skip the dependence on $(\xi,x)$ in the notation).
We assume that the derivatives of $\ti F$ with respect to $\eps$ satisfy specific bounds at $\eps=0$ (stated below in \re{assumption depsF}) and that $\tilde F(0)=0$, $\deps\ti F(0)=0$.   
$\tilde F$ will be specified later. The iteration scheme is as follows: 
The function
$(\t_0,\p_0)$, given by \re{solitonsolution}, solves
\be\la{successive eq G0}
{}&\BR{\u\dxi\l(\begin{matrix}
\t\\
\p\\
\end{matrix}\r)
-\l(\begin{matrix}
\p\\
\dx^2\t-\sin\t\\
\end{matrix}\r)
}{\large$=:{\cal G}_0(\t,\p)$}=0\,,
\ee
which is the equation characterizing the classical solitons.
In the first iteration step we amend   
${\cal G}_0(\t,\p)=0$ 
by introducing an
additional unknown variable $\lambda$ and 
adding some terms involving $(\t_0,\p_0)$ and  $\ti F$. The amended equation is of the form
\be\la{successive eq G1}
{}&\BR{\u\dxi\bma
\t\\
\p\\
\ema
-\l(\begin{matrix}
\p\\
\t\xx-\sin\t+\ti F (\eps)\\
\end{matrix}\r)
+\lambda \du\bma
\t_0\\
\p_0\\
\ema
}{\large$=:\Go(\t,\p,\lambda)$}=0\,.
\ee
Here and in the following iterations the functions $\t,\p$ depend on $(\xi,\u,x)$ and
$\lambda$ depends on $(\xi,\u)$. 
We solve ${\cal G}_1^\eps(\t,\p,\lambda)=0 $ implicitly for $(\t,\p,\lambda)$ in terms of $\eps$ and denote the solution by $(\t_1^\eps,\p_1^\eps,\lambda_{1}^\eps)$. 
%
In the next iteration step we amend ${\cal G}_1^\eps(\t,\p,\lambda)=0 $
by adding some terms involving  $(\t_1^\eps,\p_1^\eps)$ and solve the amended equation 
\be\la{successive eq G2}
{}&\BR{\u\dxi\bma
\t\\
\p\\
\ema
-\l(\begin{matrix}
\p\\
\t\xx-\sin\t+\ti F(\eps)\\
\end{matrix}\r)
+\lambda\du\bma
\t_1^0+\deps\t_1^0\eps\\
\p_1^0+\deps\p_1^0\eps\\
\ema 
}{\large$=:\Gt(\t,\p,\lambda)$}=0\,
\ee
implicitly for $(\t,\p,\lambda)$ in terms of $\eps$. 
Continuing the iteration process we obtain
in the $n$th step 
the equation
\be\la{successive eq Gn}
{}&\BR{\u\dxi\bma
\t\\
\p\\
\ema
-\l(\begin{matrix}
\p\\
\t\xx-\sin\t+\ti F(\eps)\\
\end{matrix}\r)
+\lambda\du\bma
\sum_{i=0}^{n-1} \fr{\deps^i\t_{n-1}^0}{i!}\eps^i\\
\sum_{i=0}^{n-1} \fr{\deps^i\p_{n-1}^0}{i!}\eps^i\\
\ema 
}{\large$=:\Gn(\t,\p,\lambda)$}=0\,,
\ee
where $(\t_{n-1}^\eps,\p_{n-1}^\eps,\lambda_{n-1}^\eps)$ denotes the solution of ${\cal G}_{n-1}^\eps(\t,\p,\lambda)=0 $. 
We solve $\Gn(\t,\p,\lambda)=0$ implicitly for $(\t,\p,\lambda)$ in terms of $\eps$ and 
denote the solution by $(\t_n^\eps,\p_n^\eps,\lambda_{n}^\eps)$. Due to the assumptions on $\tilde F$
it is possible to iterate this 
procedure arbitrarily many times.
The
existence of the implicit solutions $\eps\mapsto(\t_n^\eps,\p_n^\eps,\lambda_{n}^\eps)$ for $n\ge 1$
is ensured by the implicit function theorem.
In the actual proof,
we consider 
rather
the transformed equations 
\be\la{intro def tilde G}
\ti {\cal G}_n^\eps(\hatt,\hatp,\lambda):={\cal G}_n^\eps(\t_0+\hatt,\p_0+\hatp,\lambda)=0,~~~~n\ge 1,
\ee
which will be solved for $ (\hatt,\hatp,\lambda)$ in terms of $\eps$. This is caused by functional analytic reasons, among others, by the fact that $\t_0(\xi,\u,x)  \not\rightarrow 0$ as $|x|\to \infty$ for fixed $\xi$ and $\u$. 
We denote the solutions to the equations $\ti {\cal G}_n^\eps(\hatt,\hatp,\lambda)=0,~n\ge 1,$ by $( \hatt_{n}^\eps , \hatp_{n}^\eps ,\lambda_{n}^\eps)$, where $(\t_{n}^\eps ,\p_{n}^\eps ,\lambda_{n}^\eps)
 =(\t_0+\hatt_{n}^\eps ,\p_0+\hatp_{n}^\eps ,\lambda_{n}^\eps)$. 
The application of the implicit function theorem
relies on the fact that 
$(0,0,0,0)$ 
solves all equations in a particular point, i.e., 
$\ti {\cal G}_n^0(0,0,0)=0$. 
As a consequence of the construction, the 
solution obtained in the $n$th iteration
$\eps\mapsto(\t_n^\eps,\p_n^\eps,\lambda_{n}^\eps)$ solves the equation 
\be\la{intro spec PDE}
{\u\dxi\bma
\t\\
\p\\
\ema
-\l(\begin{matrix}
\p\\
\t\xx-\sin\t+\ti F (\eps)\\
\end{matrix}\r)
+\lambda \du\bma
\t\\
\p\\
\ema  
}
=0\,
\ee
up to errors of order $\eps^{n+1}$ for $ n \ge 1$.

In \cite{Mashkin}, the iterative equations
$\ti {\cal G}_n^\eps(\hatt,\hatp,\lambda)=0$
were solved in spaces of different regularity in $u$
such that the regularity of the spaces  
(which contain the corresponding iterative solutions) decreases after each iteration step by the order of $1$. This technique was used
for the following reason. 
Each iterative equation contains a derivative
with respect to $u$ of the solution of the  preceding equation, as one can see in \re{successive eq Gn}.
This derivative leads to loss of regularity in $u$ in the target set of the map $\ti {\cal G}_n$ after each iteration step.  
However,
the employment of the implicit function theorem for solving the iterative equations requires that
the corresponding linearizations are invertible and that the maps $\ti {\cal G}_n$ are well-defined.
In \cite{Mashkin}, this is ensured by considering the maps $\ti {\cal G}_n$ on spaces of decreasing regularity in $u$.
Since,
in the present paper, we need to execute infinitely many (and not only finitely many) iterations in order to obtain a sequence of implicit solutions, we 
modify the iteration scheme and proceed as follows.

Due to the analyticity assumption on $F$ in the present paper (which was not supposed in \cite{Mashkin}), the implicit solutions (as well as its derivatives) are analytic in $\eps$,
which is a consequence of the implicit function theorem.
In the first iteration we solve 
$\ti {\cal G}_1^\eps(\hatt,\hatp,\lambda)=0$
and the solution may be written in the form
\be\la{intro Taylor representation}
\begin{split}
 (\hat\t_{1 }^\eps,\hat\p_{1 }^\eps,\hat\lambda_{1 }^\eps) =  \l(\sum_{i=0}^{\infty} \fr{\deps^i\hat\t_{1 }^0}{i!}\eps^i,\sum_{i=0}^{\infty} \fr{\deps^i\hat\p_{1 }^0}{i!}\eps^i,\sum_{i=0}^{\infty} \fr{\deps^i\lambda_{1 }^0}{i!}\eps^i\r)\,
\end{split}
\ee
accordingly. 
Further application of the implicit function theorem in spaces of higher regularity in $u$
yields that $(\t_{1 }^\eps,\p_{1 }^\eps,\lambda_{1 }^\eps)$ is sufficiently often differentiable in  $u\in [-\u_*,\u_*]$, but possibly in a smaller neighbourhood of $\eps=0
$
than that where representation \re{intro Taylor representation} holds. 
We prove bounds on 
the derivatives 
$
\du^K\deps^N(\t_1^0,\p_1^0,\lambda_{1}^0)
$ (derivatives with respect to $u\in [-\u_*,\u_*]$ and $\eps $, evaluated at $\eps=0$), which have the form
\begin{align}
\la{intro first bound it1}
\forall N\ge 2,~ 0\le K \le 2:{}&&
\l\Vert\bma
\du^K\deps^N\t_1^0\\
\du^K\deps^N\p_1^0\\
\du^K\deps^N\lambda_{1}^0
\ema
\r\Vert
&\le 
 C^{2N+ 2K -3}(N-2)!,\\
\la{intro second bound it1}
\forall N\ge 2, ~K\ge 3:{}&&
\l\Vert\bma
\du^K\deps^N\t_1^0\\
\du^K\deps^N\p_1^0\\
\du^K\deps^N\lambda_{1}^0
\ema
\r\Vert
&\le 
 C^{2N + 2K-3}(N-2)!(K-3)!\,, 
\end{align}
where $\Vert \cdot \Vert$ is an appropriate norm.
These bounds 
imply that 
the implicit solution $(\t_{1 }^\eps,\p_{1 }^\eps,\lambda_{1 }^\eps)$ is differentiable in $u$ in 
the same 
neighbourhood of $\eps=0$
where also representation \re{intro Taylor representation} holds.
Thus the map $\ti {\cal G}_2$ is well defined on the same spaces where also 
$\ti {\cal G}_1^\eps(\hatt,\hatp,\lambda)=0$   
was\\ 
solved initially.
This eliminates the loss of regularity problem faced in \cite{Mashkin} (in the first iteration)
and we are able to solve
the next iterative equation 
$\ti {\cal G}_2^\eps(\hatt,\hatp,\lambda)=0$
on the same spaces as also the preceding equation $\ti {\cal G}_1^\eps(\hatt,\hatp,\lambda)=0$.
The process 
of solving the iterative equations 
will be continued using the same arguments, whereas
we prove successively
bounds on 
the derivatives of the succeeding solutions
$
\du^K\deps^N(\t_n^0,\p_n^0,\lambda_{n}^0)
$
(derivatives with respect to $u\in [-\u_*,\u_*]$ and $\eps $, evaluated at $\eps=0$).
The bounds are uniform in $n$ and have the form
\begin{align}
\retainlabel{intro first bound}
\forall N\ge 2,~ 0\le K \le 2:{}&&
\l\Vert\bma
\du^K\deps^N\t_n^0\\
\du^K\deps^N\p_n^0\\
\du^K\deps^N\lambda_{n}^0
\ema
\r\Vert
&\le 
 C^{2N+ 2K -3}(N-2)!,\\
\la{intro second bound}
\forall N\ge 2, ~K\ge 3:{}&&
\l\Vert\bma
\du^K\deps^N\t_n^0\\
\du^K\deps^N\p_n^0\\
\du^K\deps^N\lambda_{n}^0
\ema
\r\Vert
&\le 
 C^{2N + 2K-3}(N-2)!(K-3)!\,, 
\end{align}
where $\Vert \cdot \Vert$ is as above.
Here and in \re{intro first bound it1}-\re{intro second bound it1} the higher order derivatives with respect to $u$ are needed in order to control the first order derivative terms (derivative
with respect to $u$) in the iterative equations (see \re{successive eq Gn}).
This fact itself and the proof of bounds \re{intro first bound it1}-\re{intro second bound} as well rely 
on a recursive formula for 
$
\du^K\deps^N(\t_n^0,\p_n^0,\lambda_{n}^0)
$,
which is proved by induction on $N$ and $K$.
Furthermore, the assumptions on the derivatives of $\ti F$
at $\eps=0$ are used in the proof of 
\re{intro first bound it1}-\re{intro second bound}.
Bounds 
\re{intro first bound it1}-\re{intro second bound} 
imply
that
all iterative implicit solutions are defined on the same neigbourhood, 
can be represented there as Taylor series around $\eps=0$ analogous to \re{intro Taylor representation} and are there differentiable in $u$. Moreover, it follows from
\re{intro first bound it1}-\re{intro second bound} 
that the iterative implicit solutions are all contained in the same space
and that as $n\to \infty$
the sequence
$(\hat \t_n^\eps,\hat \p_n^\eps,\lambda_{n}^\eps)$ converges
to the limit
\be
(\hatt_\infty^\eps,\hatp_\infty^\eps,\lambda_\infty^\eps) := {}&\l(\sum_{i=1}^{\infty} \fr{\deps^i\t_{i}^0}{i!}\eps^i,\sum_{i=1}^{\infty} \fr{\deps^i\p_{i}^0}{i!}\eps^i,\sum_{i=0}^{\infty} \fr{\deps^i\lambda_{i}^0}{i!}\eps^i\r). 
\ee
Using these facts and \re{intro first bound it1}-\re{intro second bound} we conclude that 
the function
$$
(\t_\infty^\eps,\p_\infty^\eps,\lambda_\infty^\eps):=(\t_0+\hatt_\infty^\eps ,\p_0+\hatp_\infty^\eps ,\lambda_\infty^\eps)
$$
satisfies the equation  
\be
\la{intro eqofinterest}
\u\dxi\l(\begin{matrix}
\t_\infty^\eps\\
\p_\infty^\eps\\
\end{matrix}\r)
-\l(\begin{matrix}
\p_\infty^\eps\\
[\t_\infty^\eps]\xx-\sin\t_\infty^\eps+\ti F(\eps)\\
\end{matrix}\r)
+\lambda_{ \infty}^\eps\du\bma
\t_\infty^\eps\\
\p_\infty^\eps \\
\ema
 =0 \,.
\ee

In order to generate the invariant virtual solitary manifold, we apply the iteration scheme to a specific $\ti F$,  which is
a truncated version of the perturbation term $F$ from \re{Cauchy_intro}, given by
\be\la{intro assumption Chi and F}
\begin{cases}
 \ti F(\eps,\xi,x):=\F \chi(\xi),\\ 
 \text{where } \chi\in C^{\infty}(\R),~\chi (\xi)=1 \text{  for } |\xi|\le |\xi_s|+3 \text{ and } \chi (\xi)=0 \text{ for } |\xi|\ge |\xi_s|+4.
\end{cases}
\ee
The limit of the thereby obtained sequence of iterative solutions, defines the solution of \re{intro eqofinterest} with the specific $\tilde F$ (given by \re{intro assumption Chi and F}), which implies our main result.

In order to simplify the computations we work in the present paper on spaces, which have lower regularity in $(\xi,x)$ than the corresponding spaces in \cite{Mashkin}.

Finally let us explain under which condition
the 
parameters $\bar\xi,\bar\u$ describe
a fixed nontrivial perturbation of the uniform linear motion
as $\eps \to 0$.
We consider the setting where the 
assumption $|u_s|\le \tilde C \eps^{}$ 
is satisfied and hence where
the solution of \re{Cauchy_intro InfMF}  exists and may be expressed up to times $1/(\tilde C \eps)$ in
the mentioned way. 
For all $n\ge 1$ the linearization of 
$(\hatt,\hatp,\lambda) \mapsto 
\ti {\cal G}_n^\eps(\hatt,\hatp,\lambda)
$ 
carried out at 
$(\hatt,\hatp,\lambda)=(0,0,0)$, $\eps=0$ 
%
is invertible and we denote the linearization by
 $$
{\frak M}_0^\a:  
( \t,
\p,
\lambda) \mapsto
{\frak M}_0^\a
( \t,
\p,
\lambda ).
$$
Thus there exist functions $( \bar\t,
\bar \p,
\bar \lambda )$ such that the second derivative with respect to $\eps$ of a general function $\ti F$ (which operates on appropriate spaces), 
evaluated at $\eps=0$, can be written in the form
\be \la{intro condition on tiF}
\bma
0\\
 \deps^{2}  \ti F(0)
\ema
=
{\frak M}_0^\a
( \bar\t,
\bar \p,
\bar \lambda ),~~\text{  ${\frak M}_0^\a$ given by \cref{le invertibilityMxiCtwo alpha} (case $m=0
$)}.
\ee
Here the functions $\bar\t,\bar\p$ depend on $(\xi,\u,x)$ and $\bar\lambda$ depends on $(\xi,\u)$.
%
%
%
%
ODEs \re{intro ODEs invariant mf} can be rescaled in time by introducing $s=\eps  t
$, 
$
\hxi(s)= \bar\xi(s/\eps^{  })
$, and
$
\hu(s)= \fr 1 {\eps^{ }} {\bar u(s/\eps^{ })}
$
such that the 
corresponding transformed ODEs have the form 
\be
\ds \hxi(s) =   \hu(s) , ~~~~
\ds \hu(s) =  \fr 1 {\eps^{2 }} \lambda_{ \infty}^\eps(\hxi(s), \eps^{ }\hu(s)).\nn
\ee
As $\eps \to 0$, the transformed ODEs converge to ODEs that describe a fixed nontrivial perturbation of the uniform linear motion if 
the next condition is satisfied:
\be
\la{intro condition on nontrivial dynamic}
\begin{cases}
{}&\text{There exists $\chi$ satisfying \re{intro assumption Chi and F} such that for $\ti F$ given by \re{intro assumption Chi and F} }\\ {}& \text{the following holds: }
 \bar \lambda(\cdot, 0)\not= 0 \text{ in 
representation \re{intro condition on tiF}}.  
\end{cases}
\ee 
This is for the following reason.
The functions $(\t^\eps_\infty,\p^\eps_\infty,\lambda^\eps_\infty)$ satisfy the relation
\be\la{introduction infty relation}
\u\dxi\l(\begin{matrix}
\t_\infty^\eps\\
\p_\infty^\eps\\
\end{matrix}\r)
-\l(\begin{matrix}
\p_\infty^\eps\\
[\t_\infty^\eps]\xx-\sin\t_\infty^\eps+\ti F(\eps)\\
\end{matrix}\r)
+\lambda_{ \infty}^\eps\du\bma
\t_\infty^\eps\\
\p_\infty^\eps \\
\ema =0.
\ee
Due to the assumption on $F$ it holds that $\deps   \ti F(0)=0$
and 
differentiation of \re{introduction infty relation}
with respect to $\eps$ yields 
\be 
\bma
0\\
 \deps^l  \ti F(0)
\ema
=
{\frak M}_0^\a
( \deps^l \t_\infty^0,
\deps^l \p_\infty^0,
\deps^l\lambda_\infty^0 ),~~~~~1\le l\le 2.
\ee
Using invertibility of ${\frak M}_0^\a$,
condition \re{intro condition on nontrivial dynamic} and the fact that $\lambda_\infty^0=0$ 
it follows that $0\not=\lambda_\infty^\eps(\cdot,0)= {\cal O}(\eps^{2})$, which implies the claim.
%
%
%
\paragraph{Outline of the Paper}
The paper is organized as follows. In \cref{se: Main Results}, we formulate the main result. 
In \cref{se: Implicit function theorem}, we
modify the iteration scheme from \cite{Mashkin},
construct a sequence of iterative solutions and prove bounds on the elements of the sequence.  
In \cref{se: Convergence of the Sequence}, we show that the sequence of iterative solutions converges and that its limit satisfies the equation of interest.
Our main result, \cref{maintheorem}, is proved in \cref{se: Main Results Proof}.

\paragraph{Notation and Conventions}
For a Hilbert space $H$ we denote its inner product by $\langle\cdot,\cdot\rangle_H $.
To simplify notation,
occasionally we drop the dependence of functions on certain variables. We write ${L_x^{2}(\R)},{H_{\xi,x}^{k}(\R^2)}$ and so on for the Lebesgue and Sobolev spaces when we wish to emphasize the variables of integration.   
We use the notation $\t(\xi,\u,x)=\t(\u)(\xi,x)$, $\p(\xi,\u,x)=\p(\u)(\xi,x)$.

\newpage
\section{Main Result}\la{se: Main Results}
To formulate our result precisely, we need some definitions.
\bde\la{def:PartfourMainResult}
Let $\a,k,m\in\mathbb{N}_0$ and $u_*>0$. Let us denote by $I(\u_*):=[-\u_*,\u_*]$.
\begin{itemize}
\item [(a)] $H^{k,\alpha}(\R) $ denotes the weighted Sobolev space of functions with finite norm
$$|\t|_{H^{k,\alpha}(\R)}= |(1+|x|^2) ^\fr \alpha 2\t(x)|_{H_x^{k}(\R)}\,.$$
\item [(b)] $H^{k,\alpha}(\R^2) $ denotes the weighted Sobolev space of functions with finite norm
$$|\t|_{H^{k,\alpha}(\R^2)}= |(1+|\xi|^2+|x|^2) ^\fr \alpha 2\t(\xi,x)|_{H_{\xi,x}^{k}(\R^2)}\,.$$
\item [(c)] $ \ubar{ Y}^\a$ is the space $H^{2,\a}(\R^2) \oplus H^{1,\a}(\R^2) \oplus H^{2,\a}(\R)$
with the finite norm
$$
|y|_{\ubar{ Y}^\a} = |\t|_{H^{2,\a}(\R^2)}+ |\p|_{H^{1,\a}(\R^2)}+|\lambda |_{H^{2,\a}(\R)}\,.
$$ 
\item [(d)] $Y_m^\a(\u_*)$ is the space\\
$\\\ba
{}&\bigg\{ y=(\t,\p,\lambda) \in C^m( I(\u_*), \ubar{ Y}^\a) :  \Vert y \Vert_{Y_m^\a(\u_*)} <\infty;~\forall~ \u\in I(\u_*),~\forall~\mu\in H^{2,\a}(\R):\\
{}& \Ltwortwoaxix{ \bma \t(\u)(\xi,x)\\
\p(\u)(\xi,x) \ema}{\mu(\xi)\bma \t_K'(\g(u)(x-\xi))\\
-\u\g(u)\t_K''(\g(u)(x-\xi))\ema} =0 \bigg\}\,
\ea\\$\\
with the finite norm
$$
\Vert y \Vert_{Y_m^\a(\u_*)} =\sup_{\u\in I(\u_*)}  \l( \sum_{i=0}^m |\du^i y(u)|_{\ubar{ Y}^\a}\r) \,.
$$

\end{itemize}
\ede
\noindent The weighted Sobolev spaces in \cref{def:PartfourMainResult} (a), (b) are defined as in \cite{Kopylova}. We are now ready to state our main result.
\bth\la{maintheorem}
Let $\xi_s\in \R $, 
$\Xi:=\Xi(\xi_s):=  |\xi_s|+3$ and $\a\in \mathbb{N}_0$.
Assume that $F\in C^{\infty}((-1,1),H^{0,\a}(\R))$, $F$ is analytic and the conditions
\be \la{assumption derivatives F}
  F(0)=0,~~\deps  F(0)=0\,, 
\ee 
\be \la{assumption bounds derivatives F}
  \forall N \ge 2:~~
	\l|
	\deps^N F(0)
	\r|_{H^{0,\a}}
	\le c^N (N-2)!
	\,
\ee 
are satisfied. Then there exist $\eps^*>0$, $\u_*>0$, $\tilde C>0$ and a map 
\be\la{main theorem map}
(-\eps^*,\eps^*) \to Y_{0}^\a(\u_*),~
\eps \mapsto (\hat\t_\infty^\eps,\hat\p_\infty^\eps,\lambda_\infty^\eps)\la{map in main theorem}
\ee
of class $C^\infty$ such that the following holds. 
Let $\eps\in(0,\eps^*)$. 
Consider the Cauchy problem
\be\la{SGE1}
\partial_t \bma
\t \\
\p 
\ema
=\l(\begin{matrix}
\p \\
\dx^2\t -\sin\t +\F\\
\end{matrix}\r) ,~~~~
\bma
\t(0,x)\\
\p(0,x)
\ema=\bma
\t^\eps_\infty(\xi_s,\u_s,x)\\
\p^\eps_\infty(\xi_s,\u_s,x)
\ema ,
\ee
where $(\t_\infty^\eps,\p_\infty^\eps)=(\t_0+\hatt_\infty^\eps ,\p_0+\hatp_\infty^\eps)$ with $(\t_0,\p_0)$ given by \re{solitonsolution} such that the 
initial velocity satisfies
$|u_s|< \u_*$.
Then the Cauchy problem 
has a unique solution, which may be written in the form
\be\la{form}
\bma
\t(t,x)\\ 
\p(t,x)
\ema=\bma 
\t_\infty^\eps(\bar\xi(t),\bar\u(t),x)\\ 
\p_\infty^\eps(\bar\xi(t),\bar\u(t),x)
\ema,
\ee
where 
$\bar\xi,\bar\u$ solve the 
system of equations
\be\la{exactODE virtual1}
 {}&\d{\bar\xi}( t) =  \bar\u(t)  \,,~~
\d{\bar\u}( t) = \lambda_{\infty}^\eps\l(\bar\xi(t), \bar\u(t)\r)\,, ~~~~\bar\xi(0)=\xi_s,~\bar\u(0)=\u_s\,,
\ee
and representation \re{form} of the solution is valid
as long as $|\bar\xi(t)|\le \Xi, ~|\bar \u(t)|<\u_*$.

In particular, if $|u_s|\le \tilde C \eps^{}$,
then the Cauchy problem \re{SGE1}
has a unique solution on the time interval
\be
0\le t \le \fr {1} {\tilde C \eps } 
\ee
and may be written in the form \re{form} with ODEs \re{exactODE virtual1}. If additionally the perturbation $F$ satisfies condition \re{intro condition on nontrivial dynamic}, then  
the 
parameters $\bar\xi,\bar\u$ describe
a fixed nontrivial perturbation of the uniform linear motion 
as $\eps \to 0$. 
\eth
\noindent
The assumption on the first derivative of $F$ in \re{assumption derivatives F} is not crucial, it is made in order to simplify the computations in the proof of the bounds on the derivatives of the iterative solutions in \cref{se: Implicit function theorem} (\cref{derivatives estimate}).

We work in weighted Sobolev spaces in order  
to ensure 
that 
symplectic decomposition (implemented by techniques of \cite{Mashkin}) is possible in a neighbourhood of the invariant virtual solitary manifold, since this is promising to be useful in our future works. The well-definedness of a corresponding
symplectic orthogonality condition formulated in analogy to \cite[Theorem 2.2 (b)]{Mashkin} is guaranteed if function \re{main theorem map} maps into a weighted space $Y_{0}^\a(\u_*)$ where $\a \ge 1$ (nevertheless symplectic decomposition is not needed in the present paper).

\section{Construction of the Sequence of Iterative Solutions}\la{se: Implicit function theorem}

In this section we modify
the iteration scheme from \cite{Mashkin}
and 
construct a sequence of iterative solutions. 
By making stronger assumptions than in \cite{Mashkin} on the function $\ti F$ (utilized in the scheme below), 
we obtain more accurate information on the iterative solutions.
We start with a definition.
\bde
Let $\a,m\in\mathbb{N}_0$ and $u_*>0$.
\begin{itemize}

\item [(a)] $\ubar{ Z}^\a$ is the space $H^{1,\a}(\R^2) \oplus H^{0,\a}(\R^2)$ 
with the finite norm
$$
|z|_{\ubar{ Z}^\a} = |\v|_{H^{1,\a}(\R^2)}+ |\w|_{H^{0,\a}(\R^2)}\,.
$$

\item [(b)] $
Z_m^\a(\u_*)$ is the space $\bigg\{ z =(\v,\w) \in C^m(I(\u_*), \ubar{ Z}^\a) : \Vert z \Vert_{Z_m^\a(\u_*)} <\infty \bigg\}\,$ 
with the finite norm
$$
\Vert z \Vert_{Z_m^\a(\u_*)} =\sup_{\u\in I(\u_*)} \l( \sum_{i=0}^m |\du^i y(u)|_{\ubar{ Z}^\a}\r) \,.
$$
\item [(c)] Let us denote by
$t_{1}(\xi,\u,x):= \bma
\dxi\t_0(\xi,u,x)\\
\dxi\p_0(\xi,u,x)\\
\ema $
and by 
$t_{2}(\xi,\u,x):=\bma
\du\t_0(\xi,u,x)\\
\du\p_0(\xi,u,x)\\
\ema,
$ where $\u\in(-1,1),~\xi,x\in\R $.
\end{itemize}
\ede

\noindent
The application of the implicit function theorem in the iteration scheme is justified by the following proposition, which ensures that the corresponding linearization of 
$(\hatt,\hatp,\lambda) \mapsto 
\ti {\cal G}_n^\eps(\hatt,\hatp,\lambda),~n\ge 1,
$ 
carried out at 
$(\hatt,\hatp,\lambda)=(0,0,0)$, $\eps=0$
is invertible.

\bpro \la{le invertibilityMxiCtwo alpha}
Let $\a \in \mathbb{N}_0$.
There exists $\underline{u}^\a>0$ such that for any $m\in \mathbb{N}_0$ the operator\\ 
$
{\frak M}_m^\a:  Y_m^\a(u_*)  \to Z_m^\a(u_*),~
( \t,
\p,
\lambda) \mapsto
{\frak M}_m^\a
( \t,
\p,
\lambda ),
$
given by
\be
{\frak M}_m^\a
( \t,
\p,
  \lambda  )(\u)
=\bma
\u\dxi\t(u) -\p(u)  \\ 
-\dx^2\t(u) +\cos(\t_K(\g(u)(x-\xi)))\t(u) +\u\dxi\p(u)  \\ 
\ema
+ \lambda(u)  
t_2(\xi,\u,x),
\ee
is invertible if $0< u_*< \underline{u}^\a$.
\epro
\noindent
\bpr
The proof was given in \cite[Proposition 3.2]{Mashkin}.
\epr
\noindent
The modified iteration scheme is formalized in the following theorem. 
\bth\la{thimplicitfunctionIT1 alpha} 
Let $\a \in \mathbb{N}_0$ and let $\underline{u}^\a$
be from \cref{le invertibilityMxiCtwo alpha}. Let $0< u_*<\underline{u}^\a$,
$J=(-1,1)$ and let $\ti F: J \to H^{0,\a}(\R^2)\,, \eps \mapsto \ti F(\eps)$ be an analytic 
function such that
\be\la{assumptions tildeF}
\ti F(0)=0,~~\deps\ti F(0)=0,
\ee 
and
\be  \la{assumption depsF}
  \forall N \ge 2:~~\Bigg\Vert 
	\bma
	0\\
	\deps^N \tilde F(0)
	\ema
	\Bigg\Vert_{Z_0^\a(\u_*)} 
	\le \bar c^N (N-2)!
	\,.
\ee
\noindent
Let $\ti{\cal G}_1$ be given by
\be
{}&\ti{\cal G}_1: J  \times Y_{0}^\a (u_*) \to  Z_{0}^\a(u_*)\,,
(\eps,\hatt,\hatp,\lambda) \mapsto \tiGo(\hatt,\hatp,\lambda):=\Go(\t_0+\hatt,\p_0+\hatp,\lambda)\,,
\ee
where ${\cal G}_1$ is defined by \re{successive eq G1}. 
Then there exists $\eps^*>0$ and
a map
\be
{}&(-\eps^*,\eps^*) \to Y_{0}^\a(u_*),~
\eps \mapsto (\hat\t_1^\eps,\hat\p_1^\eps,\lambda_{ 1}^\eps)\,,
\ee
of class $C^\infty$ such that
$
\ti{\cal G}_1^\eps(\hatt_1^\eps,\hatp_1^\eps,\lambda_{1}^\eps)=0\,.
$
Let $\ti{\cal G}_2$ be given by
\be
\ti{\cal G}_2: J \times Y_{0}^\a(u_*)  \to  Z_{0}^\a(u_*)\,,
 (\eps,\hatt,\hatp,\lambda) \mapsto \tiGt(\hatt,\hatp,\lambda):=\Gt(\t_0+\hatt,\p_0+\hatp,\lambda)\,,
\ee
where ${\cal G}_2$ is defined by \re{successive eq G2} with $(\t_1^\eps ,\p_1^\eps ,\lambda_{1}^\eps)
=(\t_0+\hatt_1^\eps ,\p_0+\hatp_1^\eps ,\lambda_{1}^\eps)$. Then there exists 
a map
\be
(-\eps^*,\eps^*) \to Y_{0}^\a(u_*),
\eps \mapsto (\hat\t_2^\eps,\hat\p_2^\eps,\lambda_{ 2}^\eps)\,,
\ee
of class $C^\infty$ such that
$
\ti{\cal G}_2^\eps(\hatt_2^\eps,\hatp_2^\eps,\lambda_{2}^\eps)=0\,.
$
This process can be continued successively to arrive at $\ti{\cal G}_n$ for any 
$n\in \mathbb{N}$ 
be given by
\be
\ti{\cal G}_n: J \times Y_{0}^\a(u_*)  \to  Z_{0}^\a(u_*)\,,
(\eps,\hatt,\hatp,\lambda) \mapsto \tiGn(\hatt,\hatp,\lambda):=\Gn(\t_0+\hatt,\p_0+\hatp,\lambda)\,,
\ee 
where ${\cal G}_n$ is defined by \re{successive eq Gn} with $(\t_{n-1}^\eps ,\p_{n-1}^\eps ,\lambda_{n-1}^\eps)
=(\t_0+\hatt_{n-1}^\eps ,\p_0+\hatp_{n-1}^\eps ,\lambda_{n-1}^\eps)$.
There exists 
a map
\be
(-\eps^*,\eps^*) \to Y_{0}^\a(u_*),
\eps \mapsto (\hat\t_n^\eps,\hat\p_n^\eps,\lambda_{ n}^\eps)\,,
\ee 
of class $C^\infty$ such that
$
\ti{\cal G}_n^\eps(\hatt_n^\eps,\hatp_n^\eps,\lambda_{n}^\eps)=0\,.
$
The iterative solutions may be written in the form
\be
(\hat\t_n^\eps,\hat\p_n^\eps,\lambda_n^\eps) = {}&\l(\sum_{i=0}^{\infty} \fr{\deps^i\hat\t_{n}^0}{i!}\eps^i,\sum_{i=0}^{\infty} \fr{\deps^i\hat\p_{n}^0}{i!}\eps^i,\sum_{i=0}^{\infty} \fr{\deps^i\lambda_{n}^0}{i!}\eps^i\r)\,
\ee
as a limit in $Y_{0}^\a (u_*) $ for $\eps\in(- \eps^*, \eps^* )$. We set $(\t_{n}^\eps ,\p_{n}^\eps ,\lambda_{n}^\eps)
:=(\t_0+\hatt_{n}^\eps ,\p_0+\hatp_{n}^\eps ,\lambda_{n}^\eps)$.
\eth
\noindent
In the following we point out the relation among the derivatives of the iterative solutions from \cref{thimplicitfunctionIT1 alpha} at $\eps=0$.  
\ble \la{thITrelations}
Let the assumptions of Theorem \ref{thimplicitfunctionIT1 alpha} hold and let $n\ge 2$.\\
Then
$(\deps^k\t_{n-1}^{0},\deps^k\p_{n-1}^{0},\deps^k\lambda_{ {n-1}}^{0})=(\deps^k\t_n^{0},\deps^k\p_n^{0},\deps^k\lambda_{ n}^{0})$ for $k=0, \ldots, n-1$.
\ele
\bpr
Analogous to \cite[Theorem 3.4]{Mashkin}.
\epr
\bre
The derivatives of the iterative solutions coincide at $\eps=0$ in the following way: 
$(\deps^k\t_{1}^{0},\deps^k\p_{1}^{0},\deps^k\lambda_{ {1}}^{0})
= (\deps^k\t_{2}^{0},\deps^k\p_{2}^{0},\deps^k\lambda_{{2}}^{0})$ for $k=0, 1;$
$(\deps^k\t_{2}^{0},\deps^k\p_{2}^{0},\deps^k\lambda_{{2}}^{0})
= (\deps^k\t_{3}^{0},\deps^k\p_{3}^{0},\deps^k\lambda_{{3}}^{0})$ for $k=0, 1,2 $ and so on. 
\ere
\noindent
Now we prove some bounds on the derivatives of the iterative solutions.
These bounds will be used in the inductive proof of \cref{thimplicitfunctionIT1 alpha}.
Moreover, the bounds play a major key in the proof of convergence of the sequence of iterative solutions and they are also needed in 
order to show that the corresponding limit defines a function which satisfies
the equation of interest.

\ble\la{derivatives estimate}
Let the assumptions 
of \cref{thimplicitfunctionIT1 alpha}
be satisfied.
There exists $C>0$ such that the following holds.
Let $n\in\N$ and	
assume that for $1\le j \le n$ the iterative solutions of the equations
$\ti {\cal G}_j^\eps(\hatt,\hatp,\lambda)=0$
%
exist, then 
the following bounds are satisfied:

\begin{align}
\la{lebound1}
1\le K \le 2:{}&&
\l\Vert\bma
\du^K\t_0\\
\du^K\p_0\\
0
\ema
\r\Vert_{Y_0^\a(\u_*)} 
&\le 
C ,
\\
\la{lebound2}
\forall K\ge 3:{}&&
\l\Vert\bma
\du^K\t_0\\
\du^K\p_0\\
0
\ema
\r\Vert_{Y_0^\a(\u_*)}
&\le 
C^{2K-3}(K-3)!,
\\
\la{lebound3}
\forall N\ge 2,~ 0\le K \le 2:{}&&
\l\Vert\bma
\du^K\deps^N\t_n^0\\
\du^K\deps^N\p_n^0\\
\du^K\deps^N\lambda_{n}^0
\ema
\r\Vert_{Y_0^\a(\u_*)} 
&\le 
 C^{2N+ 2K -3}(N-2)!,
\\
\la{lebound4}
\forall N\ge 2, ~K\ge 3:{}&&
\l\Vert\bma
\du^K\deps^N\t_n^0\\
\du^K\deps^N\p_n^0\\
\du^K\deps^N\lambda_{n}^0
\ema
\r\Vert_{Y_0^\a(\u_*)} 
&\le 
 C^{2N + 2K-3}(N-2)!(K-3)!\,.
%
\end{align}
\ele
\bpr 
An argument for differentiability with respect to $u$ of the iterative solutions will be given in the proof of \cref{thimplicitfunctionIT1 alpha}. 
The upper bounds in this proof are given by sums of certain types and the 
major key is that those sums converge.  
In the following we take a closer look at one of them, since the other cases can be treated similarly.
It holds for $l\ge 6$ that
\be
{}&
\sum_{k=3}^{l-3}\fr{(l-1) (l-2)}{(l-1-k)!k!} 
(k-3)! (l-k-3)!
\\
={}&
\sum_{k=3}^{l-3}\fr{(l-1) (l-2)}{(l-1-k)(l-2-k)k(k-1)(k-2)} 
 \\
={}&
\sum_{3 \le k \le \lfloor (l-1)/2 \rfloor } \fr{(l-1) (l-2)}{(l-1-k)(l-2-k)k(k-1)(k-2)} 
\\
{}&
+\sum_{ \lfloor (l-1)/2 \rfloor< k \le l-3} 
\fr{(l-1) (l-2)}{(l-1-k)(l-2-k)k(k-1)(k-2)} 
\ee
\be
\le{}&
\sum_{3 \le k \le \lfloor (l-1)/2 \rfloor } 
 \fr{ 1}{\fr{(l-1-k)}{(l-1)}\fr{(l-2-k)}{(l-2)}k(k-1)(k-2)} 
%
\\
{}&
+
\sum_{   2\le j < l-1 -\lfloor (l-1)/2 \rfloor} 
\fr{1}{j(j-1)\fr{(l-1-j)}{l-1}\fr{(l-2-j)}{(l-2)}(l-3-j)} 
\\
\le{}&
\sum_{3 \le k \le \lfloor (l-1)/2 \rfloor } 
\fr{4}{ (k-2)(k-1)k}
+
\sum_{   2\le j < l-1 -\lfloor (l-1)/2 \rfloor} 
\fr{4}{j(j-1) }=:R(l)\,
\ee
and thus $\sup_l R(l)<\infty $.
Let us now deduce a recursive relation which will be needed later. 
Taking the 
$K$-th derivative with respect  to $u$ of ${\cal G}_0(\t_0,\p_0)=0$ yields
\be
0{}&=\bma
\u\dxi\du^K\t_0-\du^K\p_0\\
\u\dxi\du^K\p_0-\dx^2\du^K\t_0
\ema
+\bma
0\\
\sum_{m=1}^{K-1}\binom{K-1}{m}  \du^{m} \cos(\t_0) \du^{K-m}\t_0 +  \cos(\t_0) \du^{K}\t_0
\ema\\
{}&
+K\bma
\dxi\du^{K-1}\t_0\\
\dxi\du^{K-1}\p_0
\ema
\,.
\ee
Thus
\be\la{rekursivrelationKderivative}
\du^K\bma
 \t_0\\
 \p_0\\
0
\ema
=
-  \l[{\frak M}_{0}^\a\r]^{-1} \Bigg[
\bma
0\\
\sum_{m=1}^{K-1}\binom{K-1}{m}  \du^{m} \cos(\t_0) \du^{K-m}\t_0 
\ema
+K\bma
\dxi\du^{K-1}\t_0\\
\dxi\du^{K-1}\p_0
\ema
\Bigg]
\,.
\ee
\noindent
We show first \re{lebound1}-\re{lebound2}.
We chose $C>1$ such that the claim \re{lebound1}-\re{lebound2} is true for $0\le K \le 3$ and such that  
$
\sup_{u\in I(u_*)} |\du^{m}\cos\t_0|_{L^\infty_{\xi,x}(\R^2)}  \le 
C \, 
$
for $0\le m \le 3$.
In the following we will put some more assumptions on $C$, where we tag each of them with an exclamation mark "!". 
We assume that the claim \re{lebound1}-\re{lebound2} holds for all integers up to 
$K-1$ and prove the induction step. Let $n\in \mathbb{N}$.
Firstly, we show that for $3\le m \le K$:
\be\la{ind cos} 
\sup_{u\in I(u_*)}|\du^{m}\cos\t_0|_{L^\infty_{\xi,x}(\R^2)}  \le 
(m-3)!C^{2m-3+1/3}\,. 
\ee
We assume that \re{ind cos} holds for all integers $3\le m \le K-1$ and show the
induction step. In the following we use Sobolev embedding theorems.
Notice that
\be
{}&
\sup_{u\in I(u_*)}  \Bigg|
\sum_{k=0}^{l-1}\binom{l-1}{k} \du^k\cos(\t_0) \du^{l-k}\ten 
\Bigg|_{L^\infty_{\xi,x}(\R^2)} 
\\
={}&\sup_{u\in I(u_*)}  \Bigg| \Bigg(\sum_{k=3}^{l-3}\fr{(l-1)!}{(l-1-k)!k!} \du^k\cos(\t_0) \du^{l-k}\t_0
+\cos(\t_0)\du^{l}\t_0
+(l-1)\du\cos(\t_0)\du^{l-1}\t_0
\\
{}&
+\fr{(l-1)(l-2)}2\du^2\cos(\t_0)\du^{l-2}\t_0
+\du^{l-1}\cos(\t_0)\du\t_0
+(l-1)\du^{l-2}\cos(\t_0)\du^{2}\t_0
\Bigg)
\Bigg|_{L^\infty_{\xi,x}(\R^2)}\\
\ee
\be
\le{}&
(l-3)!\sum_{k=3}^{l-2}\fr{(l-1)(l-2)}{(l-1-k)!k!}  (k-3)!(l-k-3)! 
C^{2k-3+1/3} C^{2(l-k)-3}\\
{}&+
(l-3)!C^{2l-3}
+(l-1)C^{}
(l-5)!C^{2(l-1)-3}
+3\fr{(l-1)(l-2)}2   C^{4-3+1/3} 
(l-5)!C^{2(l-2)-3}
\\
{}&
+
(l-4)!C^{2(l-3)}C+
(l-1)(l-5)!C^{2(l-2)-3+1/3}C\\
\OT{\le}{ ! }{}&
(l-3)!C^{2l-3+1/3} \,.
\ee
Using this estimate it follows for $3\le m \le K$ that
\be
{}&
\sup_{u\in I(u_*)}  \Bigg|
\du^m(\cos(\t_0))
\Bigg|_{L^\infty_{\xi,x}(\R^2)} 
\\
={}&
\sup_{u\in I(u_*)}  \Bigg|
\du^{m-1}(\sin(\t_0)\deps\t_0)
\Bigg|_{L^\infty_{\xi,x}(\R^2)} 
\\
={}&
\sup_{u\in I(u_*)}  \Bigg|
\sum_{l=0}^{m-1}\binom{m-1}{l}  \du^{l} \sin(\t_0) \du^{m-l}\t_0 
\Bigg|_{L^\infty_{\xi,x}(\R^2)} 
\\
= {}&
\sup_{u\in I(u_*)}  \Bigg|
\Bigg(\sum_{l=1}^{m-1}\binom{m-1}{l}  \du^{l-1} \l(\cos(\t_0) \du\t_0\r)\du^{m-l}\t_0 + \sin\t_0\du^m\t_0\Bigg)
\Bigg|_{L^\infty_{\xi,x}(\R^2)} 
\\
= {}&
\sup_{u\in I(u_*)}  \Bigg|
\Bigg(\sum_{l=1}^{m-1}\binom{m-1}{l} \l( \sum_{k=0}^{l-1}\binom{l-1}{k} \du^k\cos(\t_0) \du^{l-k}\t_0\r)\du^{m-l}\t_0 + \sin\t_0\du^m\t_0\Bigg)
\Bigg|_{L^\infty_{\xi,x}(\R^2)} 
\\
= {}&
\sup_{u\in I(u_*)}  \Bigg|
\Bigg(\sum_{l=3}^{m-1}\binom{m-1}{l} \l( \sum_{k=0}^{l-1}\binom{l-1}{k} \du^k\cos(\t_0) \du^{l-k}\t_0\r)\du^{m-l}\t_0 \\
{}&+(m-1) \cos(\t_0)\du\t_0
+ \sin\t_0\du^m\t_0\\
{}&+\fr{(m-1)(m-2)}2 \cos(\t_0)\du^2\t_0+\fr{(m-1)(m-2)}2\du\cos(\t_0)\du\t_0\Bigg)
\Bigg|_{L^\infty_{\xi,x}(\R^2)} 
\\
\le {}& (m-3)!\sum_{l=3}^{m-1}\fr{(m-1)(m-2)}{(m-l-1)!l!} (l-3)! (m-l-3)!
C^{2l-3+1/3}C^{2(m-l)-3}\\
{}&+ (m-1)C+ 
(m-2)!
C^{2m-3}+
\fr{(m-1)(m-2)}2 
C^3
+\fr{(m-1)(m-2)}2CC\\
\OT{\le}{ ! }{}& 
(m-3)! C^{2m-3+1/3}  \,,
\ee
which completes the induction step for \re{ind cos}. 
In the following we denote by $\Vert \cdot \Vert$ the operator norm of $\l[{\frak M}_{0}^\a\r]^{-1} $.
Now we estimate 
$
\du^K(
 \t_0,
 \p_0,
 0)
$
by using the recursive formula \re{rekursivrelationKderivative} and the bounds 
\re{ind cos}:
\be
{}&\l\Vert
\l[{\frak M}_{0}^\a\r]^{-1} 
\l[
\bma
0\\
 \sum_{m=1}^{K-1}\binom{K-1}{m}  \du^{m} \cos(\t_n^0) \du^{K-m}\t_n^0
\ema
+K\bma
\dxi\du^{K-1}\t_0\\
\dxi\du^{K-1}\p_0
\ema\r]
\r\Vert_{Y_0^\a(\u_*)}\\
\le{}&
\l\Vert\l[{\frak M}_{0}^\a\r]^{-1}  \r\Vert
\Bigg((K-3)!
\sum_{m=3}^{K-1}\fr{(K-1)(K-2)}{(K-m-1)!m!} (m-3)!(K-m-3)!
C^{2m-3+1/3} C^{2(K-m)-3}\\
{}&
+(K-1)(K-4)!
CC^{2(K-1)-3}
+
\fr{(K-1)(K-2)(K-5)!}2 C^{4-3+1/3}C^{2(K-2)-3}\\
{}&
+
 (K-4)!CC^{2(K-1)-3+1/3 }
+
 (K-1) (K-5)!CC^{2(K-2)-3+1/3  } 
+K(K-4)!C^{2(K-1)-3  } 
\Bigg)\\
\OT{\le}{!}{}&
(K-3)!
C^{2K-3-1/3}\,.
\ee
Assuming that $C^{2K-3-1/3}\OT{\le}{ ! } C^{2K-3}$, the induction step for \re{lebound1}-\re{lebound2} is complete.
\noindent
Before proving the remaining claim, we deduce some recursive relations for further computations.
Taking the 
$N$-th derivative with respect  to $\eps$ of $\Gn(\ten,\pen,\luen)=0$ 
yields
\be\la{ITn}
0={}&\deps^N\Gn(\ten,\pen,\luen)\\
={}&
\bma
\u\dxi\deps^N\ten-\deps^N\pen\\
\u\dxi\deps^N\pen-\dx^2\deps^N\ten
\ema
+\bma
0\\
\sum_{m=1}^{N-1}\binom{N-1}{m}  \deps^{m} \cos(\ten) \deps^{N-m}\ten +  \cos(\ten) \deps^{N}\ten 
\ema\\
{}&
-\bma
0\\
 \deps^N\tiF 
\ema
+ \bma
\sum_{i=0}^{n-1} \sum_{l=0}^{N}\binom{N}{l}  \deps^{N-l} \luen \deps^{l} \l[\fr{\du\deps^i\t_{n-1}^0}{i!}\eps^i\r]\\
\sum_{i=0}^{n-1} \sum_{l=0}^{N}\binom{N}{l}  \deps^{N-l} \luen \deps^{l} \l[\fr{\du\deps^i\p_{n-1}^0}{i!}\eps^i\r]\\
\ema\,.
\ee
\noindent
Thus we obtain 
\be\la{recursive relation N}
\bma
\deps^N\t_n^0\\
\deps^N\p_n^0\\
\deps^N\lambda_{n}^0
\ema
={}&\l[{\frak M}_{0}^\a\r]^{-1} \Bigg[\bma
0\\
\deps^N \ti F(0)
\ema
-
\bma
 \sum_{1\le l \le \min\{n-1,N-1\}} \binom{N}{l}  \deps^{N-l} \lambda_{n}^0 {\du\deps^l\t_{n}^0}\\
 \sum_{1\le l \le \min\{n-1,N-1\}} \binom{N}{l}  \deps^{N-l} \lambda_{n}^0 {\du\deps^l\p_{n}^0}
\ema
\\
{}&
-\bma
0\\
 \sum_{m=1}^{N-1}\binom{N-1}{m}  \deps^{m} \cos(\ten) \deps^{N-m}\ten
\ema
\Bigg|_{\eps=0}
\Bigg]
\,.
\ee
Due to assumption \re{assumptions tildeF}
it follows from case $N=1$ combined with \cref{le invertibilityMxiCtwo alpha} 
that 
$
(\deps\t_n^0,
\deps\p_n^0,
\deps\lambda_{n}^0)=(0,0,0)$.
Taking the $K$th derivative with respect to $u$
of \re{ITn} yields\\
\be
0={}&
\bma
\u\dxi\du^K\deps^N\t_n^0-\du^K\deps^N\p_n^0\\
\u\dxi\du^K\deps^N\p_n^0-\dx^2\du^K\deps^N\t_n^0
\ema
+
\bma
0\\
\cos(\t_n^0)\du^{K}\deps^{N}\t_n^0 
\ema
+\du^{K}  \deps^{N} \lambda_{n}^0
\bma
\du \t_0\\
\du\p_0
\ema
\\
{}&
+
\sum_{\substack{0\le m \le N-1,\\ 0 \le k \le K, ~(m,k)\not=(0,0) }}\binom{N-1}{m} \binom{K}{k}
\bma
0\\
  \du^{k} \deps^{m} \cos(\ten) \du^{K-k}\deps^{N-m}\ten 
\ema
\Bigg|_{\eps=0}
\\
{}&
+ \sum_{\substack{0\le l \le \min\{n-1,N\}\\
0\le k \le K, ~(l,k)\not=(0,0) }} \binom{N}{l}\binom{K}{k}  \bma
  \du^{K-k}  \deps^{N-l} \lambda_{n}^0 {\du^{k+1} \deps^l\t_{n}^0}\\
  \du^{K-k}  \deps^{N-l} \lambda_{n}^0 {\du^{k+1} \deps^l\p_{n}^0}
\ema
+K\bma
\dxi\du^{K-1}\deps^N\t_n^0\\
\dxi\du^{K-1}\deps^N\p_n^0
\ema
\, .
\ee
Thus 
we obtain
\be\la{recursive relation KN}
\begin{split}
{}&\bma
\du^K\deps^N\t_n^0\\
\du^K\deps^N\p_n^0\\
\du^K\deps^N\lambda_{n}^0
\ema
\\
={}&-\l[{\frak M}_{0}^\a\r]^{-1} \Bigg[
\sum_{\substack{0\le m \le N-1,\\ 0 \le k \le K, ~(m,k)\not=(0,0) }}\binom{N-1}{m}  \binom{K}{k}
\bma
0\\
  \du^{k} \deps^{m} \cos(\ten) \du^{K-k}\deps^{N-m}\ten 
\ema
\Bigg|_{\eps=0}
\\
{}&
+\sum_{\substack{0\le l \le \min\{n-1,{ N-1}\}\\
0\le k \le K, ~(l,k)\not=(0,0) }} 
\binom{N}{l}\binom{K}{k}  \bma
  \du^{K-k}  \deps^{N-l} \lambda_{n}^0 {\du^{k+1} \deps^l\t_{n}^0}\\
  \du^{K-k}  \deps^{N-l} \lambda_{n}^0 {\du^{k+1} \deps^l\p_{n}^0}
\ema
+K\bma
\dxi\du^{K-1}\deps^N\t_n^0\\
\dxi\du^{K-1}\deps^N\p_n^0
\ema
\Bigg]
.
\end{split}
\ee
Now we show \re{lebound3}-\re{lebound4}.
We prove the claim by induction on $N$, whereas we conduct for each $N$ an induction on $K$. 
In some further estimates we will use the fact that there exists $c>0$ such that 
$$
 |\lambda\t|_{H^{1,\alpha}(\R^2)} 
\le  
c|\lambda|_{H^{2,\alpha}(\R)}| \t|_{H^{1,\alpha}(\R^2)}\, 
$$
for $\lambda\in H^{2,\alpha}(\R)$ and $\t \in H^{1,\alpha}(\R^2) $. This follows from Morrey's inequality. Let us start the induction.

\noindent\underline{$N=1$:} The terms 
$(\du^K\deps\t_n^0,
\du^K\deps\p_n^0,
\du^K\deps\lambda_{n}^0)$ vanish for any $K$ due to assumption \re{assumptions tildeF}. 
\\
\noindent\underline{$N=2$:} This case can be treated similarly to the following proof of the induction step.
\noindent
\underline{$2,\ldots,N-1\rightarrow N$:}
We assume that 
bound \re{lebound3} holds for derivatives with respect to $\eps$ of order $2$ up to order $N-1$ and for derivatives with respect to $u$ of order $0$ up to order $2$.
Moreover, we assume that 
bound \re{lebound4} holds for derivatives with respect to $\eps$ of order $2$ up to order $N-1$ and for all derivatives with respect to $u$  from order $3$.
Now we show the
induction step $2,\ldots,N-1\rightarrow N$. This will be done by induction on $K$, where we use \re{recursive relation N} and \re{recursive relation KN}.\\
\underline{$K=0$}:
We consider separately the terms of 
the recursive formula \re{recursive relation N}. 
Due to \re{assumption depsF} we are able to estimate\\
\be 
{}&\l\Vert
\l[{\frak M}_{0}^\a\r]^{-1} 
\l[
\bma
0\\
 \deps^N \ti F(0)
\ema
\r]
\r\Vert_{Y_0^\a(\u_*)}\\
\OT{\le}{ ! }{}&
(N-2)!
C^{2N-3-1/3}\,,
\\
{}&\l\Vert
\l[{\frak M}_{0}^\a\r]^{-1} 
\l[
\bma
0\\
 \sum_{m=1}^{N-1}\binom{N-1}{m}  \deps^{m} \cos(\ten) \deps^{N-m}\ten
\ema
\r]
\Bigg|_{\eps=0}\r\Vert_{Y_0^\a(\u_*)}\\
\le{}&
\l\Vert\l[{\frak M}_{0}^\a\r]^{-1}  \r\Vert
\Bigg((N-2)!
\sum_{m=3}^{N-2}\fr{(N-1)}{(N-m-1)!m!} (m-2)!(N-m-2)!C^{2m-3} C^{2(N-m)-3}\\
{}&
+\fr{(N-1)(N-2)(N-4)!}2 C^{2\cdot 2-3+}C^{2(N-2)-3}
\Bigg)
\\
\OT{\le}{ ! }{}&
(N-2)!
C^{2N-3-1/3}\,,
\\
{}&
\l\Vert
\l[{\frak M}_{0}^\a\r]^{-1} 
\bma
 \sum_{l=1}^{N-1}\binom{N}{l}  \deps^{N-l} \lambda_{n-1}^0 {\du\deps^l\t_{n-1}^0}\\
 \sum_{l=1}^{N-1}\binom{N}{l}  \deps^{N-l} \lambda_{n-1}^0 {\du\deps^l\p_{n-1}^0}
\ema
\r\Vert_{Y_0^\a(\u_*)}\\
\le{}& 
\l\Vert
\l[{\frak M}_{0}^\a\r]^{-1} 
\r\Vert
(N-2)!
\sum_{l=2}^{N-2}\fr{N(N-1)}{(N-l)!l!} 
(N-l-2)! (l-2)!C^{2(N-l)-3}C^{2l-3}\\
\OT{\le}{ ! }{}&
(N-2)!
C^{2N-3-1/3}\,.
\ee
Further we assume that $3C^{2N-3-1/3}\OT{\le}{ ! } C^{2N-3}$.
\\
\underline{$K=1$}:
We consider separately the terms of 
the recursive formula \re{recursive relation KN} and obtain 
\be
{}&\Bigg\Vert
\l[{\frak M}_{0}^\a\r]^{-1} 
\Bigg[
\sum_{\substack{0\le m \le N-1,\\ 0 \le k \le 1, ~(m,k)\not=(0,0) }}\binom{N-1}{m}  \binom{1}{k}
\bma
0\\
  \du^{k} \deps^{m} \cos(\ten) \du^{1-k}\deps^{N-m}\ten 
\ema
\Bigg]
\Bigg|_{\eps=0}\Bigg\Vert_{Y_0^\a(\u_*)}\\
\le{}&
\l\Vert\l[{\frak M}_{0}^\a\r]^{-1} \r\Vert
\Bigg((N-2)!
\sum_{\substack{2\le m \le N-2,\\ 0 \le k \le 1, ~(m,k)\not=(0,0) }} \binom{1}{k} \fr{(N-1)}{(N-m-1)!m!} (m-2)!(N-m-2)!C^{2N-6} \\
{}&
  + (N-2)! C C^{2N-3} 
\Bigg)\\
\OT{\le}{ ! }{}&
(N-2)!  C^{2N+2\cdot 1-3-1/3} \,,
\ee

\be
{}&\Bigg\Vert
\l[{\frak M}_{0}^\a\r]^{-1} 
\Bigg[
\sum_{\substack{0\le l \le \min\{n-1,N-1\}\\
0\le k \le 1, ~(l,k)\not=(0,0) }} 
\binom{N}{l}\binom{1}{k}  \bma
  \du^{1-k}  \deps^{N-l} \lambda_{n}^0 {\du^{k+1} \deps^l\t_{n}^0}\\
  \du^{1-k}  \deps^{N-l} \lambda_{n}^0 {\du^{k+1} \deps^l\p_{n}^0}
\ema
+\bma
\dxi\deps^N\t_n^0\\
\dxi\deps^N\p_n^0
\ema
\Bigg]
\Bigg\Vert_{Y_0^\a(\u_*)}\\
\le{}&
\l\Vert\l[{\frak M}_{0}^\a\r]^{-1}  \r\Vert
\Bigg((N-2)!
\sum_{\substack{2\le m \le N-2,\\ 0 \le k \le K, ~(m,k)\not=(0,0) }} \binom{K}{k} \fr{N(N-1)}{(N-m)!m!} (m-2)!(N-m-2)!C^{2N-6} \\
{}& 
+ (N-2)! C C^{2N-3} 
+(N-2)!C^{2N-3}
\Bigg)\\
\OT{\le}{ ! }{}&
(N-2)!C^{2N+2\cdot 1-3-1/3}\,.
\ee
Further we assume that $2C^{2N+2-3-1/3}\OT{\le}{ ! } C^{2N+2-3}$.
\\
\underline{$K=2$}:
We consider separately the terms of 
the recursive formula \re{recursive relation KN} and obtain
\be
{}&\Bigg\Vert
\l[{\frak M}_{0}^\a\r]^{-1} 
\Bigg[
\sum_{\substack{0\le m \le N-1,\\ 0 \le k \le 2, ~(m,k)\not=(0,0) }}\binom{N-1}{m}  \binom{2}{k}
\bma
0\\
  \du^{k} \deps^{m} \cos(\ten) \du^{2-k}\deps^{N-m}\ten
\ema
\Bigg]
\Bigg|_{\eps=0}\Bigg\Vert_{Y_0^\a(\u_*)}\\
\le{}&
\l\Vert\l[{\frak M}_{0}^\a\r]^{-1}  \r\Vert
\Bigg((N-2)!
\sum_{\substack{2\le m \le N-2,\\ 0 \le k \le 2, ~(m,k)\not=(0,0) }} \binom{2}{k} \fr{(N-1)}{(N-m-1)!m!} (m-2)!(N-m-2)!C^{2N-6} \\
{}&
  + (N-2)!  2 C C^{2N+2\cdot 1-3} + (N-2)! C C^{2N-3} 
\Bigg)\\
\OT{\le}{ ! }{}& (N-2)!  C^{2N+2\cdot 2-3-1/3} \,,
\ee
\be
{}&\Bigg\Vert
\l[{\frak M}_{0}^\a\r]^{-1} 
\Bigg[
\sum_{\substack{0\le l \le \min\{n-1,N-1\}\\
0\le k \le 2, ~(l,k)\not=(0,0) }} 
\binom{N}{l}\binom{2}{k}  \bma
  \du^{2-k}  \deps^{N-l} \lambda_{n}^0 {\du^{k+1} \deps^l\t_{n}^0}\\
  \du^{2-k}  \deps^{N-l} \lambda_{n}^0 {\du^{k+1} \deps^l\p_{n}^0}
\ema
+2\bma
\dxi\du^{ }\deps^N\t_n^0\\
\dxi\du^{ }\deps^N\p_n^0
\ema
\Bigg]
\Bigg\Vert_{Y_0^\a(\u_*)}\\
\le{}&
\l\Vert\l[{\frak M}_{0}^\a\r]^{-1}  \r\Vert
\Bigg((N-2)!
\sum_{\substack{2\le m \le N-2,\\ 0 \le k \le 2, ~(m,k)\not=(0,0) }} \binom{2}{k} \fr{N(N-1)}{(N-m)!m!} (m-2)!(N-m-2)!C^{2N-6} \\
{}& 
+  2 (N-2)! C C^{2N+2\cdot 1-3} + (N-2)! C C^{2N-3}+ 2 (N-2)!C^{2N+2\cdot 1-3}
\Bigg)\\
\OT{\le}{ ! }{}& (N-2)!C^{2N+2\cdot 2-3-1/3}.
\ee
Further we assume that $2C^{2N+2-3-1/3}\OT{\le}{ ! } C^{2N+2-3}$.
\\
\underline{$K=3$}: 
This case can be proven analogously to the case $K=2$.
 
\noindent
\underline{$0,\ldots,K-1\rightarrow K$}: 
We assume that the claim holds for all integers up to $K-1$ and show the
induction step.
Recall that in the case $N=0$ we have proven:
$$
0\le k\le 2: ~\sup_{u\in I(u_*)}| \du^{k}\cos\t_0
|_{L^\infty_{\xi,x}(\R^2)} 
 \le 
C\,,~~~~
\forall k\ge 3: ~
\sup_{u\in I(u_*)}| 
\du^{k}\cos\t_0 
|_{L^\infty_{\xi,x}(\R^2)} 
\le 
(k-3)!C^{2k-3+1/3}\,.
$$
To begin with, we show that for $2\le m \le N-1$:
\be \la{ind km cos one} 
0\le k\le 2:~~ \sup_{u\in I(u_*)}| \du^{k}\deps^{m}\cos\ten  |_{\eps=0}
|_{L^\infty_{\xi,x}(\R^2)} 
\le
 (m-2)!C^{2m+ 2k-3+1/3}\,, 
\ee
\be\la{ind km cos two} 
\forall k\ge 3:~~
\sup_{u\in I(u_*)}|
\du^{k}\deps^{m}\cos\ten   |_{\eps=0}
|_{L^\infty_{\xi,x}(\R^2)} 
\le 
(k-3)!(m-2)!C^{2k+2m-3+1/3}\,. 
\ee
The induction basis for $N=2$ can be shown similarly to the case $N=0$.
We assume that \re{ind km cos one}-\re{ind km cos two} holds for all integers $2\le m \le N-2$ and show the
induction step. 
%
We start with a preliminary estimate for $l\ge 4,~i\ge 3$:
\be
{}&
\sup_{u\in I(u_*)}\Bigg|
\sum_{k=0}^{l-1}\binom{l-1}{k} \sum_{j=0}^{i}\binom{i}{j}\du^j\deps^k\cos(\ten) \du^{i-j}\deps^{l-k}\ten
\Bigg|_{\eps=0}
\Bigg|_{L^\infty_{\xi,x}(\R^2)} 
\\
=
{}&
\sup_{u\in I(u_*)}\Bigg|
\Bigg(\sum_{k=2}^{l-2}\binom{l-1}{k} \sum_{j=3}^{i-3}\binom{i}{j}\du^j\deps^k\cos(\ten) \du^{i-j}\deps^{l-k}\ten\\
{}&+\cos(\ten) \du^{i}\deps^{l}\ten 
+i\du\cos(\ten) \du^{i-1}\deps^{l}\ten \\
{}&+\fr{i(i-1)}{2}\du^2\cos(\ten) \du^{i-2}\deps^{l}\ten 
+\sum_{j=3}^{i}\binom{i}{j}\du^j\cos(\ten) \du^{i-j}\deps^{l}\ten\\
{}&
+\sum_{k=2}^{l-2}\binom{l-1}{k} \deps^k\cos(\ten) \du^{i}\deps^{l-k}\ten
\\
{}&
+\sum_{k=2}^{l-2}\binom{l-1}{k} i\du\deps^k\cos(\ten) \du^{i-1}\deps^{l-k}\ten
+\sum_{k=2}^{l-2}\binom{l-1}{k} \fr{i(i-1)}{2}\du^2\deps^k\cos(\ten) \du^{i-2}\deps^{l-k}\ten\\
{}&+\sum_{k=2}^{l-2}\binom{l-1}{k} \fr{i(i-1)}{2}\du^{i-2}\deps^k\cos(\ten) \du^{2}\deps^{l-k}\ten
+\sum_{k=2}^{l-2}\binom{l-1}{k}i\du^{i-1}\deps^k\cos(\ten) \du^{}\deps^{l-k}\ten\\
{}&+\sum_{k=2}^{l-2}\binom{l-1}{k} \du^i\deps^k\cos(\ten) \deps^{l-k}\ten
\Bigg)\Bigg|_{\eps=0}\Bigg|_{L^\infty_{\xi,x}(\R^2)}  
\ee
\be
\le{}&
(l-2)!(i-3)!C^{2i+2l-6+2/3}\cdot\\
{}&\sum_{k=3}^{l-2}
\sum_{j=3}^{i-3}
\fr{(l-1)i(i-1)(i-2)}{(l-1-k)k(k-1)(i-j)(i-j-1)(i-j-2)j(j-1)(j-2)}   \\
%
{}&+(i-3)!(l-2)! C^{2i+2l-3}
+i (i-4)!(l-2)! CC^{2(i-1)+2l-3}  \\
{}&+\fr{i(i-1)(i-5)!(l-2)!}{2}C C^{2(i-2)+2l-3}\\
{}&
+(l-2)!(i-3)!\sum_{j=3}^{i}\fr{i(i-1)(i-2)}{(i-j)(i-j-1)(i-j-2)j(j-1)(j-2)} C^{2i+2l-6}\\
{}&
%
+(l-2)!(i-3)!\sum_{k=2}^{l-2} \fr{l-1}{k(k-1)(l-k-1)}C^{2i+2l-6}
\\
{}&
+(l-2)!i(i-4)!\sum_{k=2}^{l-2} \fr{(l-1)}{k(k-1) (l-k-1)} C^{2i+2l-6}
\\
{}&
+(l-2)!\fr{i(i-1)}{2} (i-5)!\sum_{k=2}^{l-2} \fr{(l-1)}{k(k-1) (l-k-1)} C^{2i+2l-6}
\\
{}&
+(l-2)!\fr{i(i-1)}{2} (i-5)!\sum_{k=2}^{l-2} \fr{(l-1)}{k(k-1) (l-k-1)} C^{2i+2l-6}
\\
{}&
+(l-2)!i(i-4)!\sum_{k=2}^{l-2} \fr{(l-1)}{k(k-1) (l-k-1)} C^{2i+2l-6}
\\
%
{}&
+(l-2)!(i-3)!\sum_{k=2}^{l-2} \fr{l-1}{k(k-1)(l-k-1)}C^{2i+2l-6}\\
\OT{\le}{ ! } {}&(l-2)!(i-3)!C^{2i+2l-3} \,.\la{preliminaryboundcos}
\ee
\noindent
Applying the Leibniz's formula we obtain for $2\le m \le N-1$, $ r\ge 3$: 
\be
{}&
\sup_{u\in I(u_*)}\Bigg|
\du^r\deps^m(\cos(\ten))\Big|_{\eps=0}
\Bigg|_{L^\infty_{\xi,x}(\R^2)} 
\\
={}&
\sup_{u\in I(u_*)}\Bigg|
\du^r\deps^{m-1}(\sin(\ten)\deps\ten)\Big|_{\eps=0}
\Bigg|_{L^\infty_{\xi,x}(\R^2)} 
\ee
\be
={}&
\sup_{u\in I(u_*)}\Bigg|
\sum_{l=0}^{m-1}\binom{m-1}{l}   \sum_{i=0}^{r}\binom{r}{i}\du^i\deps^{l} \sin(\ten) \du^{r-i}\deps^{m-l}\ten \Bigg|_{\eps=0}
\Bigg|_{L^\infty_{\xi,x}(\R^2)} 
\\
= {}&
\sup_{u\in I(u_*)}\Bigg|
\Bigg(\sum_{l=1}^{m-1}\binom{m-1}{l}
\sum_{i=0}^{r}\binom{r}{i}
\du^i\deps^{l-1} \l(\cos(\ten) \deps\ten\r)\du^{r-i}\deps^{m-l}\ten\\
{}&
 + \sum_{i=0}^{r}\binom{r}{i}\du^i\sin\ten\du^{r-i}\deps^m\ten\Bigg)\Bigg|_{\eps=0}
\Bigg|_{L^\infty_{\xi,x}(\R^2)} 
\\
= {}&
\sup_{u\in I(u_*)}\Bigg|
\Bigg(\sum_{l=1}^{m-1}\binom{m-1}{l}\sum_{i=0}^{r}\binom{r}{i} \l[ \sum_{k=0}^{l-1}\binom{l-1}{k} \sum_{j=0}^{i}\binom{i}{j} \du^j\deps^k\cos(\ten) \du^{i-j}\deps^{l-k}\ten\r]\du^{r-i}\deps^{m-l}\ten\\
{}&
+ \sum_{i=0}^{r}\binom{r}{i}\du^i\sin\ten\du^{r-i}\deps^m\ten\Bigg)\Bigg|_{\eps=0}
\Bigg|_{L^\infty_{\xi,x}(\R^2)} \,.
\ee
In order to control the expression
\be\la{costerm1}
\sum_{l=1}^{m-1}\sum_{i=0}^{r}\binom{m-1}{l}\binom{r}{i} \l[ \sum_{k=0}^{l-1}\binom{l-1}{k} \sum_{j=0}^{i}\binom{i}{j} \du^j\deps^k\cos(\ten) \du^{i-j}\deps^{l-k}\ten\r]\du^{r-i}\deps^{m-l}\ten
\Bigg|_{\eps=0}
\ee  
we split the sum \re{costerm1} over indices $l,i$ into two sums, one over indices $I_{m,r}:=\{(l,i) ~:~ 3\le l \le m-2  ~\text{and}~  3\le i \le r-2\}$ and the other over indices $ \{(l,i) ~:~ 1\le l \le m-1,~ 0\le i \le r\} \setminus I_{m,r}$. Using bound \re{preliminaryboundcos} for the square brackets term we estimate the sum over indices $I_{m,r}$ 
by
\be
{}& 
\sup_{u\in I(u_*)}\Bigg|
\sum_{l=3}^{m-2}\binom{m-1}{l}\sum_{i=3}^{r-2}\binom{r}{i} \l( \sum_{k=0}^{l-1}\binom{l-1}{k} \sum_{j=0}^{i}\binom{i}{j} \du^j \deps^k\cos(\ten)  \du^{i-j}\deps^{l-k}\ten\r)
\cdot\\
{}& 
\du^{r-i}\deps^{m-l}\ten
\Bigg|_{\eps=0}
  \Bigg|_{L^\infty_{\xi,x}(\R^2)} \\
\le
{}&
(m-2)!(r-3)!C^{2r+2m-9}\cdot\\
{}& 
\sum_{l=3}^{m-2}\sum_{i=3}^{r-2}\fr{(m-1) r(r-1)(r-2)}{(m-l-1)!l! (r-i)!i!} (l-2)! (i-3)! (r-i-3)!(m-l-2)!
\\
\le 
{}&
(m-2)!(r-3)! C^{2r+2m-9} 
\cdot 
\\
{}& 
    \sum_{l=3}^{m-2}\sum_{i=3}^{r-2}
\fr{(m-1) r(r-1)(r-2)}{(m-l-1) l(l-1) (r-i)(r-i-1) (r-i-2)i(i-1)(i-2)},   
\ee
where the supremum over $(r,m)$ of the double sum is finite. 
All terms of the sum over indices 
\be
{}&
\{(l,i) ~:~ 1\le l \le m-1,~ 0\le i \le r\} \setminus I_{m,r}\la{setofindices}\\
={}&
\{(l,0),~l= 1,\ldots,m-1\} \cup
\{(l,1),~l= 1,\ldots,m-1 \} \cup
\{(l,2),~l= 1,\ldots,m-1 \} \\
{}& 
\cup 
\{(l,r-1),~l= 1,\ldots,m-1 \} \cup
\{(l,r-2),~l= 1,\ldots,m-1 \} \cup
\{(1,i),~i= 0,\ldots,r\} \\
{}& 
\cup 
\{(2,i),~i= 3,\ldots,r\} \cup
\{(m-1,i),~i= 3,\ldots,r\}\cup
\{(m-2,i),~i= 3,\ldots,r\}
\ee
can be treated in a similar way, whereby one considers separately the sums over the subsets above. For instance, for indices $\{(l,0),~l= 1,\ldots,m-1\} $, we obtain due to \re{assumption derivatives F}
\be
{}& 
\sup_{u\in I(u_*)}\Bigg|
\sum_{l=1}^{m-1}\binom{m-1}{l}
\sum_{k=0}^{l-1}\binom{l-1}{k}     \deps^k\cos(\ten)  \deps^{l-k}\ten
%
 %
\du^{r }\deps^{m-l}\ten
\Bigg|_{\eps=0}
\Bigg|_{L^\infty_{\xi,x}(\R^2)}
\\
\le
{}& 
\sup_{u\in I(u_*)}\Bigg|
\binom{m-1}{2}  \sum_{k=0}^{1}\binom{1}{k}     \deps^k\cos(\ten)  \deps^{2-k}\ten
%
 %
\du^{r }\deps^{m-2}\ten
\Bigg|_{\eps=0}
\Bigg|_{L^\infty_{\xi,x}(\R^2)}\\
{}&+
\sup_{u\in I(u_*)}\Bigg|
\binom{m-1}{3}  \sum_{k=0}^{2}\binom{2}{k}     \deps^k\cos(\ten)  \deps^{3-k}\ten
%
 %
\du^{r }\deps^{m-3}\ten
\Bigg|_{\eps=0}
\Bigg|_{L^\infty_{\xi,x}(\R^2)}\\
{}&
+\sup_{u\in I(u_*)}\Bigg|
\sum_{l=4}^{m-1}\binom{m-1}{l} 
\sum_{k=0}^{l-1}\binom{l-1}{k}     \deps^k\cos(\ten)  \deps^{l-k}\ten
%
 %
\du^{r }\deps^{m-l}\ten
\Bigg|_{\eps=0}
\Bigg|_{L^\infty_{\xi,x}(\R^2)}\\
\le
{}& 
\sup_{u\in I(u_*)}\Bigg|
\fr{(m-1)(m-2)}{2}    \cos(\ten)  \deps^{2 }\ten
%
 %
\du^{r }\deps^{m-2}\ten
\Bigg|_{\eps=0}
\Bigg|_{L^\infty_{\xi,x}(\R^2)}\\
{}&+
\sup_{u\in I(u_*)}\Bigg|
\fr{(m-1)(m-2)(m-3)}{6} 
%
\cos(\ten)  \deps^{3}\ten
%
%
 %
\du^{r }\deps^{m-3}\ten
\Bigg|_{\eps=0}
\Bigg|_{L^\infty_{\xi,x}(\R^2)}\\
{}&
+\sup_{u\in I(u_*)}\Bigg|
\sum_{l=4}^{m-1}\binom{m-1}{l} \Bigg( \sum_{k=2}^{l-1}\binom{l-1}{k}     \deps^k\cos(\ten)  \deps^{l-k}\ten
+
 \cos(\ten)\deps^{l}\ten
\Bigg)
\du^{r }\deps^{m-l}\ten
\Bigg|_{\eps=0}
\Bigg|_{L^\infty_{\xi,x}(\R^2)}\\
\le
{}& 
\fr{(m-1)(m-2)}{2} (m-4)!      (r-3)!C^{1+2r+2(m-2)-3 }\\
%
{}&+
\fr{(m-1)(m-2)(m-3)}{6} 
%
(m-5)! (r-3)!C^{1+2r+2(m-3)-3 }\\
{}& +(m-2)!(r-3)! C^{2r+2m-8} \cdot 
\\
{}& 
  \Bigg( 
	%
%
+ \sum_{l=4}^{m-2}\sum_{k=2}^{l-1}
\fr{(m-1)}{(m-l-1) (l-k-1) k(k-1)}
+  \sum_{l=3}^{m-2}
\fr{(m-1)}{(m-l-1) l(l-1)}\Bigg) ,
\ee   
where the supremum over $(m,l)$ of the expression in the last line is finite. \\
Now we consider the sum 
\be
%
 \sum_{i=0}^{r}\binom{r}{i}\du^i\sin\ten\du^{r-i}\deps^m \ten
\Bigg|_{\eps=0}.
%
\ee
In order to control this sum, we write it 
by utilizing Leibniz's formula in the following way:
\be
{}&
\Bigg(\sin\ten\du^r\deps^m\ten+
 \sum_{i=5}^{r}\binom{r}{i} 
\sum_{p=0}^{i-1} \binom{i-1}{p}\du^{p}(\cos\ten) \du^{i-p}\ten
\\
{}&
+\sum_{i=1}^{4}\binom{r}{i} 
\sum_{p=0}^{i-1} \binom{i-1}{p}\du^{p}(\cos\ten) \du^{i-p}\ten\la{terms1}
\\
{}&
+\sum_{i=r-2}^{r}\binom{r}{i} 
\sum_{p=0}^{i-1} \binom{i-1}{p}\du^{p}(\cos\ten) \du^{i-p}\ten \Bigg)
\Bigg|_{\eps=0}\la{terms2}
.
\ee
Using the induction hypothesis we estimate the first term by 
$$
\sup_{u\in I(u_*)} |
\sin\ten\du^r\deps^m\ten
\big|_{\eps=0}
|_{L^\infty_{\xi,x}(\R^2)}
\le
(m-2)!(r-3)! C^{2r+2m-3} .
$$
For the second term we obtain 
\be
{}&
\sup_{u\in I(u_*)}\Bigg| 
 \sum_{i=5}^{r-3}\binom{r}{i} 
\sum_{p=0}^{i-1} \binom{i-1}{p}\du^{p}(\cos\ten) \du^{i-p}\ten
%
\du^{r-i}\deps^m\ten
\Bigg|_{\eps=0}
\Bigg|_{L^\infty_{\xi,x}(\R^2)}\\
={}&
\sup_{u\in I(u_*)}\Bigg|
 \sum_{i=5}^{r-3}\binom{r}{i} \Bigg( \sum_{p=3}^{i-3} \binom{i-1}{p}\du^{p}(\cos\ten) \du^{i-p}\ten
+ (\cos\ten) \du^{i}\ten
+ (i-1)\du^{ }(\cos\ten) \du^{i-1}\ten
\\
{}& 
+ \fr{(i-1)(i-2)}2\du^{2}(\cos\ten) \du^{i-2}\ten
+(i-1) \du^{i-2}(\cos\ten) \du^{2}\ten
\\
{}& 
+ \du^{i-1}(\cos\ten) \du^{}\ten
\Bigg)\du^{r-i}\deps^m\ten
\Bigg|_{\eps=0}
\Bigg|_{L^\infty_{\xi,x}(\R^2)}
\\
\le
{}&
(m-2)!(r-3)!C^{2r+2m-8} \cdot \\
{}& \sum_{i=5}^{r-3}\Bigg(\sum_{p=3}^{i-3}\fr{r(r-1)(r-2) }{(r-i)(r-i-1)(r-i-2)i (i-p-1)(i-p-2)p(p-1)(p-2)}   
\\
{}&
+
\fr{r(r-1)(r-2)}{i(i-1)(i-2)(r-i) (r-i-1)(r-i-2)}
+
\fr{r(r-1)(r-2)}{i(i-2)(i-3)(r-i) (r-i-1)(r-i-2)}
%
\ee
\be
{}&
+
\fr{r(r-1)(r-2)}{2i (i-3)(i-4)(r-i) (r-i-1)(r-i-2)}
\\
{}&
+
\fr{r(r-1)(r-2)}{i (i-2)(i-3)(i-4)(r-i) (r-i-1)(r-i-2)} \\
{}&
+
\fr{r(r-1)(r-2)}{i (i-1)(i-2)(i-3)(r-i) (r-i-1)(r-i-2)}
\Bigg), 
\ee
where the supremum of the sum over $r$ is finite. 
The summands of the sums \re{terms1}-\re{terms2} can be treated similarly. As an example we consider the case $i=2$: 
\be
{}&
\sup_{u\in I(u_*)}\Bigg|
\binom{r}{2}\sum_{p=0}^{1} \du^{p}(\cos\ten) \du^{2-p}\ten
\du^2\sin\ten\du^{r-2}\deps^m\ten
\Bigg|_{\eps=0}
\Bigg|_{L^\infty_{\xi,x}(\R^2)}
\\
\le 
{}&
(m-2)! r(r-1) (r-5)! C^{2r+2m-8}.
\ee
This completes the induction step for \re{ind km cos two}, since due to previous estimates an appropriate constant $C$ can be found as it was done in the cases $0\le K \le 2$. One shows \re{ind km cos one} similarly.
\noindent
Now we estimate separately the terms of 
the recursive formula \re{recursive relation KN}. 
Firstly, we start for $K\ge 5, ~N\ge 3$ with the term
\be\la{firsttermnorm}
 \l\Vert
\l[{\frak M}_{0}^\a\r]^{-1} 
\Bigg[
\sum_{\substack{0\le m \le N-1,\\ 0 \le k \le K, ~(m,k)\not=(0,0) }}\binom{N-1}{m}  \binom{K}{k}
\bma
0\\
  \du^{k} \deps^{m} \cos(\ten) \du^{K-k}\deps^{N-m}\ten 
\ema\Bigg]
\Bigg|_{\eps=0}\r\Vert_{Y_0^\a(\u_*)}.
\ee
We split the sum over indices $m,k$, analogous to \re{costerm1}, into two sums, one over indices $J_{N,K}:=\{(m,k) ~:~ 3\le m \le N-1  ~\text{and}~  3\le k \le K\}$ and the other over indices $ \{(m,k)\not=0 ~:~ 0\le m \le N-1,~ 0\le k \le K\} \setminus J_{N,K}$.
The sum over indices $J_{N,K}$ can be estimated by 
\be
{}&\l\Vert
\l[{\frak M}_{0}^\a\r]^{-1} 
\Bigg[
\sum_{\substack{3\le m \le N-1,\\ 3 \le k \le K, ~
}}\binom{N-1}{m}  \binom{K}{k}
\bma
0\\
  \du^{k} \deps^{m} \cos(\ten) \du^{K-k}\deps^{N-m}\ten
\ema\Bigg]
\Bigg|_{\eps=0}\r\Vert_{Y_0^\a(\u_*)}\\
\le {}&\l\Vert\l[{\frak M}_{0}^\a\r]^{-1}  \r\Vert (N-2)!(K-3)!C^{2K+2N-5}\cdot\\
{}&\sum_{\substack{3\le m \le N-1,\\ 3 \le k \le K 
}}\fr{(N-1)K(K-1)(K-2)}{(N-m-1)!m!(K-k)!k!} (k-3)!(m-2)!(K-k-3)!(N-m-2)!
\ee
\be
\le 
{}&\l\Vert\l[{\frak M}_{0}^\a\r]^{-1}  \r\Vert
(N-2)!(K-3)!C^{2K+2N-5}\cdot
\\
{}&
\sum_{\substack{3\le m \le N-1,\\ 3 \le k \le K
 }}\fr{(N-1)}{(N-m-1)m(m-1)}
\fr{K(K-1)(K-2)}{(K-k)(K-k-1)(K-k-2)k(k-1)(k-2)}\\
\OT{\le}{!}{}&
(N-2)!(K-3)!C^{2K+2N-4}.
\ee
We decompose the set of indices 
$ \{(m,k)\not=0 ~:~ 0\le m \le N-1,~ 0\le k \le K\} \setminus J_{N,K}$
analogously to \re{setofindices} and consider the sums over the corresponding subsets. All those sums can be treated similarly. For instance, for indices $\{(2,k),~k= 0,\ldots,K\} $, we obtain by using \re{ind km cos one}-\re{ind km cos two}: 
\be
{}&
\Bigg\Vert
 \l[{\frak M}_{0}^\a\r]^{-1}   
 \fr{(N-1)(N-2)}2 \sum_{k=0}^K \binom{K}{k}
\bma
0\\
  \du^{k} \deps^{ 2} \cos(\ten) \du^{K-k}\deps^{N-2}\ten
\ema 
\Bigg|_{\eps=0}
\Bigg\Vert_{Y_0^\a(\u_*)}\\
\le 
{}&
\l\Vert\l[{\frak M}_{0}^\a\r]^{-1} \r\Vert \Bigg(
%
 \fr{(N-1)(N-2)}2 \Bigg\Vert \sum_{k=3}^K \binom{K}{k}
\bma
0\\
  \du^{k} \deps^{ 2} \cos(\ten) \du^{K-k}\deps^{N-2}\ten 
\ema
\Bigg|_{\eps=0}
\Bigg\Vert_{Z_0^\a(\u_*)}
\\
{}&+\fr{(N-1)(N-2)}2  (K-3)!(N-4)! C^{2+2K+2(N-2)-3}\\
 {}&+\fr{(N-1)(N-2)}2 K (K-4)!(N-4)! C^{4+2(K-1)+2(N-2)-3} \\
 {}& +  \fr{(N-1)(N-2)}2  (K-5)!(N-4)! C^{6+2(K-2)+2(N-2)-3}   
\Bigg) \\
\OT{\le}{!}{}&
(N-2)!(K-3)!C^{2K+2N-4}.
\ee
Secondly, we consider for $K\ge 5, ~N\ge 3$ the term
\be
{}&\l\Vert
\l[{\frak M}_{0}^\a\r]^{-1} 
\Bigg[
\sum_{\substack{0\le l \le \min\{n-1,N-1\}\\
0\le k \le K, ~(l,k)\not=(0,0) }} 
\binom{N}{l}\binom{K}{k}  \bma
  \du^{K-k}  \deps^{N-l} \lambda_{n}^0 {\du^{k+1} \deps^l\t_{n}^0}\\
  \du^{K-k}  \deps^{N-l} \lambda_{n}^0 {\du^{k+1} \deps^l\p_{n}^0}
\ema
  \Bigg]
\r\Vert_{Y_0^\a(\u_*)} ,
\ee
from \re{recursive relation KN}. We treat this term analogously to \re{firsttermnorm} and
the sum over indices $J_{N,K}$ can be estimated by
\be
{}&
\l\Vert
\l[{\frak M}_{0}^\a\r]^{-1} 
\Bigg[
\sum_{\substack{3\le l \le \min\{n-1,N-1\}\\
3\le k \le K
}} 
\binom{N}{l}\binom{K}{k}  \bma
  \du^{K-k}  \deps^{N-l} \lambda_{n}^0 {\du^{k+1} \deps^l\t_{n}^0}\\
  \du^{K-k}  \deps^{N-l} \lambda_{n}^0 {\du^{k+1} \deps^l\p_{n}^0}
\ema
  \Bigg]
\r\Vert_{Y_0^\a(\u_*)}
%
\ee
\be
\le
{}&\l\Vert\l[{\frak M}_{0}^\a\r]^{-1}  \r\Vert (N-2)!(K-3)! C^{2K+2N-5} \cdot\\
{}&\sum_{\substack{3\le m \le N-1,\\ 3 \le k \le K 
}}\fr{N(N-1)K(K-1)(K-2)}{(N-m)!m!(K-k)!k!} (k-2)!(m-2)!(K-k-3)!(N-m-2)!
\\
\le
{}&\l\Vert\l[{\frak M}_{0}^\a\r]^{-1}  \r\Vert (N-2)!(K-3)!  C^{2K+2N-5} \cdot \\
{}&\sum_{\substack{3\le m \le N-1,\\ 3 \le k \le K
 }}\fr{N(N-1)}{(N-m)(N-m-1)m(m-1)}
\fr{K(K-1)(K-2)}{(K-k)(K-k-1)(K-k-2)k(k-1)}\\
\OT{\le}{!}{}&
(N-2)!(K-3)!C^{2K+2N-4}\,.
 \ee
We decompose the set of indices 
$ \{(m,k)\not=0 ~:~ 0\le m \le N-1,~ 0\le k \le K\} \setminus J_{N,K}$ and estimate the corresponding sums as above. For instance, for indices $\{(0,k),~k= 1,\ldots,K\} $, we obtain 
\be
{}&
\Bigg\Vert
 \l[{\frak M}_{0}^\a\r]^{-1}   
\sum_{k=1}^K \binom{K}{k}
\bma
\du^{K-k}  \deps^{N} \lambda_{n}^0 {\du^{k+1} \t_{n}^0}\\
  \du^{K-k}  \deps^{N} \lambda_{n}^0 {\du^{k+1} \p_{n}^0}
\ema  
\Bigg\Vert_{Y_0^\a(\u_*)}\\
\le 
{}&
\l\Vert\l[{\frak M}_{0}^\a\r]^{-1} \r\Vert \Bigg( \Bigg\Vert
\sum_{k=3}^K \binom{K}{k}
\bma
\du^{K-k}  \deps^{N} \lambda_{n}^0 {\du^{k+1} \t_{n}^0}\\
  \du^{K-k}  \deps^{N} \lambda_{n}^0 {\du^{k+1} \p_{n}^0}
\ema 
\Bigg\Vert_{Z_0^\a(\u_*)} \\
{}&
+  K (K-4)!(N-2)!C^{2(K-1)+2N-5}  + K (K-1) (K-5)!(N-2)! C^{2(K-2)+2N-3}
\Bigg) \\
\OT{\le}{!}{}& 
(N-2)!(K-3)!C^{2K+2N-3-1/3} .
\ee
The last term in \re{recursive relation KN}, 
\be
{}&\l\Vert
\l[{\frak M}_{0}^\a\r]^{-1} 
\Bigg[
K\bma
\dxi\du^{K-1}\deps^N\t_n^0\\
\dxi\du^{K-1}\deps^N\p_n^0
\ema
  \Bigg]
\r\Vert_{Y_0^\a(\u_*)}\,,
\ee
can be estimated by
\be
{}&
\l\Vert\l[{\frak M}_{0}^\a\r]^{-1} \r\Vert
\Bigg\Vert
K\bma
\dxi\du^{K-1}\deps^N\t_n^0\\
\dxi\du^{K-1}\deps^N\p_n^0
\ema
\Bigg\Vert_{Z_0^\a(\u_*)}\\
\le{}&
\l\Vert\l[{\frak M}_{0}^\a\r]^{-1} \r\Vert
 K(K-4)!(N-2)!C^{2(K-1)+2N-3}\\
\OT{\le}{!}{}&
(N-2)!(K-3)!C^{2K+2N-4}\,,
\ee
which completes the proof by the same argument as in the cases $0\le K \le 2$.
\epr
\noindent
By using \cref{derivatives estimate} we prove now \cref{thimplicitfunctionIT1 alpha}.
\bpr[of \cref{thimplicitfunctionIT1 alpha} ]
In this proof, we use the notation $Y_m^\a=Y_m^\a(\u_*)$, $Z_m^\a=
Z_m^\a(\u_*)$.
We refer to \cite[Theorem 15.1]{Deimling} and check their proof of the implicit function theorem, whereas we show that $r$ and ${\delta}$ do not depend on $\ti{\cal G}_n$.
Once
$
\ti{\cal G}_{n }: J \times Y_{0}^\a  \to  Z_{0}^\a
$ 
is defined, one obtains that its derivative 
with respect to $(\hatt,\hatp,\lambda)$ evaluated at $(\eps,\hatt,\hatp,\lambda)=(0,0,0,0)$ 
is given by
%
%
%
${\frak M}_{0}^\a$. We set 
$$
S_{n }(\eps,\hatt,\hatp,\lambda)=  \l[{\frak M}_{0}^\a\r]^{-1}  \ti{\cal G}_{n}^\eps(\hatt,\hatp,\lambda)-I(\hatt,\hatp,\lambda)\,.
$$
We start with $\ti{\cal G}_{1}$. Notice that
$
\ti{\cal G}_1^{0}(0,0,0)
=0\,.
$  
Let the constant $C$ be such that it satisfies the assumptions demanded in the proof of \cref{derivatives estimate}.
Since 
$
D_{(\hatt,\hatp,\lambda)}S_{1 }(0,0,0,0)=0 
$
and $D_{(\hatt,\hatp,\lambda)}S_{1 }$ is continuous, we fix $k\in(0,1)$ and find $1\ge\delta>0$ such that
\be\la{bound DS1}
\l\Vert 
D_{(\hatt,\hatp,\lambda)}S_{1 }(\eps,\hatt,\hatp,\lambda)
\r\Vert_{Z_0^\a(\u_*)} + \l\Vert \l[{\frak M}_{0}^\a\r]^{-1}\r\Vert \sum_{n=1}^\infty c_n\eps^n
\le k
\ee
on $\overline B_{\delta }(0) \times \overline  B_{\delta }(0) $, where
$
c_1= C$,
$c_n= \fr{ C^{2n-3}}{n(n-1)}
$ 
for $n\ge 2$
and 
$\Vert \cdot \Vert$ denotes the operator norm of $\l[{\frak M}_{0}^\a\r]^{-1} $.
Since 
$
S_{1 }(0,0,0,0)=0 
$
and 
$
S_{1 }(\cdot,0,0,0)
$
is continuous, there exists $r=:\bar\eps\le\delta$ such that 
$$
\l\Vert 
S_{1 }(\eps,0,0,0)
\r\Vert_{Z_0^\a(\u_*)}
< \delta(1-k)
$$
on 
$\overline B_{r}(0)$.
Thus there exists by \cite[Theorem 15.1]{Deimling}
a map
$$
(-\bar\eps,\bar\eps) \to Y_{0}^\a,
~~~\eps \mapsto (\hat\t_{1}^\eps,\hat\p_{1}^\eps,\lambda_{1}^\eps)
$$
such that
$
\tiG_{1}(\hat\t_{1}^\eps,\hat\p_{1}^\eps,\lambda_{1}^\eps)=0\,. 
$
Let $\underaccent{\bar} \eps>0$ be the radius of convergence of $\sum_{n=2}^\infty c_n\eps^n$ and $\eps^* := \min \{  \underaccent{\bar} \eps, \bar \eps\}$. 
Since $\ti F$ is analytic, the solution $(\hat\t_{1}^\eps,\hat\p_{1}^\eps,\lambda_{1}^\eps)$ is also analytic and 
may be written in the form
\be\la{sol_zero}
(\hat\t_{1 }^\eps,\hat\p_{1 }^\eps,\lambda_{1 }^\eps) = {}&\l(\sum_{i=0}^{\infty} \fr{\deps^i\hat\t_{1 }^0}{i!}\eps^i,\sum_{i=0}^{\infty} \fr{\deps^i\hat\p_{1 }^0}{i!}\eps^i,\sum_{i=0}^{\infty} \fr{\deps^i\lambda_{1 }^0}{i!}\eps^i\r)\,
\ee
for $\eps\in(-\eps^*,\eps^* )$
due to \cref{derivatives estimate}.	
%
Considering the map
$\ti{\cal G}_{1,m}$ on spaces of higher regularity, given by
\be
\ti{\cal G}_{1,m}: J \times Y_{m}^\a   \to  Z_{m}^\a \,,
(\eps,\hatt,\hatp,\lambda) \mapsto \tiGo(\hatt,\hatp,\lambda):=\Go(\t_0+\hatt,\p_0+\hatp,\lambda)\,,
\ee 
where ${\cal G}_1$ is defined by \re{successive eq G1},
we
obtain in the same way for any $m\in {\mathbb N}$ a constant
$\bar\eps_{m}>0$ and a map  
$$
(-\bar\eps_{m},\bar\eps_{m}) \to Y_{m}^\a,
~~~\eps \mapsto (\hat\t_{1,m}^\eps,\hat\p_{1,m}^\eps,\lambda_{1,m}^\eps)
$$
such that
$
\tiG_{1,m}(\hat\t_{1,m}^\eps,\hat\p_{1,m}^\eps,\lambda_{1,m}^\eps)=0\,.
$
Since 
$\ti F$ is analytic and
$
(
\hat\t_{1 }^\eps,
 \hat\p_{1 }^\eps,
 \lambda_{1 }^\eps
)=(\hat\t_{1,m}^\eps,\hat\p_{1,m}^\eps,\lambda_{1,m}^\eps)\in Y_{m}^\a \,
$ 
for $\eps\in(-\bar\eps_{m},\bar\eps_{m}) $, it follows from \cref{derivatives estimate} that 
$
(
\hat\t_{1 }^\eps,
 \hat\p_{1 }^\eps,
 \lambda_{1 }^\eps
)\in Y_{m}^\a \,
$ for $\eps\in(- \eps^*, \eps^* )$ and consequently that $
\ti{\cal G}_{2 }: J \times Y_{0}^\a  \to  Z_{0}^\a
$ 
is well defined. 
Since 
\be
0=\l\Vert\bma
\deps^1\t_1^0\\
\deps^1\p_1^0\\
\deps^1\lambda_{1}^0
\ema
\r\Vert_{Y_0^\a(\u_*)}\le c_1
\ee
due to \re{assumption derivatives F},\re{recursive relation N} and \cref{le invertibilityMxiCtwo alpha},
it follows from \re{bound DS1} that
$$
\l\Vert 
D_{(\hatt,\hatp,\lambda)}S_{2 }(\eps,\hatt,\hatp,\lambda)
\r\Vert_{Z_0^\a(\u_*)} 
\le k
$$
on $\overline B_{\delta }(0) \times \overline  B_{\delta }(0) $. 
Obviously
$$
\l\Vert 
S_{2 }(\eps,0,0,0)
\r\Vert_{Z_0^\a(\u_*)}
< \delta(1-k)
$$
on 
$\overline B_{r}(0)$.
Thus there exists by the same argument as above
an analytic map
$$
(-\eps_{ }^*,\eps_{ }^*) \to Y_{0}^\a,
~~~\eps \mapsto (\hat\t_{2 }^\eps,\hat\p_{2 }^\eps,\lambda_{2 }^\eps),
$$
which may be written in a form analogous to \re{sol_zero} for $\eps\in(- \eps^*,  \eps^* )$ such that  
$
\tiG_{2 }(\hat\t_{2 }^\eps,\hat\p_{2 }^\eps,\lambda_{2 }^\eps)=0\,.
$
We continue this process successively, whereas we use in the second and in the succeeding iteration steps the following argument. 
Assuming that the first $n-1$ iterative solutions are obtained, it holds
\be
\fr 1 {N!}\l\Vert\bma
\deps^N\t_{n-1}^0\\
\deps^N\p_{n-1}^0\\
\deps^N\lambda_{n-1}^0
\ema
\r\Vert_{Y_0^\a(\u_*)}\le c_{N}~~~~\text{for}~~~~1\le N\le n-1,
\ee
due to 
\cref{derivatives estimate}. 
Thus \re{bound DS1} yields that
$$
\l\Vert 
D_{(\hatt,\hatp,\lambda)}S_{n }(\eps,\hatt,\hatp,\lambda)
\r\Vert_{Z_0^\a(\u_*)} 
\le k
$$
on $\overline B_{\delta }(0) \times \overline  B_{\delta }(0) $. 
Since
$$
\l\Vert 
S_{n }(\eps,0,0,0)
\r\Vert_{Z_0^\a(\u_*)}
< \delta(1-k)
$$
on 
$\overline B_{r}(0)$
there exists by the same argument as above
an analytic map
\be
{}&(-\eps_{ }^*,\eps_{ }^*) \to Y_{0}^\a,
~~~\eps \mapsto (\hat\t_{n }^\eps,\hat\p_{n }^\eps,\lambda_{n}^\eps),
\ee 
which may be written in a form analogous to \re{sol_zero} for $\eps\in(-  \eps^*,  \eps^* )$ such that
$
\tiG_{n }(\hat\t_{n }^\eps,\hat\p_{n }^\eps,\lambda_{n }^\eps)=0.
$ 
\epr
\section{Convergence of the Sequence of Iterative Solutions}\la{se: Convergence of the Sequence}
In this section, we show that the sequence of iterative solutions constructed in \cref{se: Implicit function theorem} converges and that its limit defines a function which solves the equation of interest.
\ble\la{le limit}
Let $\a$, $u_*$ and $\eps^*$ be from \cref{thimplicitfunctionIT1 alpha}. The limit
\be
(\hatt_\infty^\eps,\hatp_\infty^\eps,\lambda_\infty^\eps) := {}&\l(\sum_{i=1}^{\infty} \fr{\deps^i\t_{i}^0}{i!}\eps^i,\sum_{i=1}^{\infty} \fr{\deps^i\p_{i}^0}{i!}\eps^i,\sum_{i=0}^{\infty} \fr{\deps^i\lambda_{i}^0}{i!}\eps^i\r) 
\ee
exists in $Y_{0}^\a (u_*) $ for $\eps\in(-\eps^*,\eps^* )$. 
We set $(\t_\infty^\eps,\p_\infty^\eps,\lambda_\infty^\eps):=(\t_0+\hatt_\infty^\eps ,\p_0+\hatp_\infty^\eps ,\lambda_\infty^\eps)$ with $(\t_0,\p_0)$ given by \re{solitonsolution}.
\ele
\bpr
The claim follows from \cref{thimplicitfunctionIT1 alpha} and \cref{derivatives estimate}, since $\eps^*$ is less or equal than the radius of convergence of
$
\sum_{n=2}^\infty \fr{C^{2n-3}}{n(n-1)} \eps^n
$ with $C$ from \cref{derivatives estimate}.
\epr

\bth\la{le limitsolution}
Let $u_*$ and $\eps^*$ be from \cref{thimplicitfunctionIT1 alpha}.
Then it holds for any $ \u\in I(u_*) $ and $\eps\in(-\eps^*, \eps^* )$ that     
$$\ba\la{ITnquantitativ}
{}&
\u\dxi\l(\begin{matrix}
\t_\infty^\eps\\
\p_\infty^\eps\\
\end{matrix}\r)
-\l(\begin{matrix}
\p_\infty^\eps\\
[\t_\infty^\eps]\xx-\sin\t_\infty^\eps+\ti F(\eps)\\
\end{matrix}\r)
+\lambda_{ \infty}^\eps\du\bma
\t_\infty^\eps\\
\p_\infty^\eps \\
\ema =0\\
\ea.
$$
\eth
\bpr
Let $n\in  \mathbb{N}$.
Notice that $$\forall \u\in I(u_*):~~
\u\dxi\l(\begin{matrix}
\t_n^\eps\\
\p_n^\eps\\
\end{matrix}\r)
-\l(\begin{matrix}
\p_n^\eps\\
[\t_n^\eps]\xx-\sin\t_n^\eps+\ti F(\eps)\\
\end{matrix}\r)
+\lambda_{ n}^\eps\du\bma
\sum_{i=0}^{n-1} \fr{\deps^i\t_{n}^0}{i!}\eps^i\\
\sum_{i=0}^{n-1} \fr{\deps^i\p_{n}^0}{i!}\eps^i \\
\ema =0\,.
$$
It holds due to \cref{thimplicitfunctionIT1 alpha} that
\be
(\hat\t_n^\eps,\hat\p_n^\eps,\lambda_n^\eps) = {}&\l(\sum_{i=0}^{\infty} \fr{\deps^i\hat\t_{n}^0}{i!}\eps^i,\sum_{i=0}^{\infty} \fr{\deps^i\hat\p_{n}^0}{i!}\eps^i,\sum_{i=0}^{\infty} \fr{\deps^i\lambda_{n}^0}{i!}\eps^i\r)\,.
\ee
Thus using \cref{thITrelations} and \cref{derivatives estimate} we obtain for $n\ge 2$
and $\eps\in(-\eps^*, \eps^* )$:
\be
{}&\l\Vert \bma
\t_\infty^\eps -\t_n^\eps\\
\p_\infty^\eps -\p_n^\eps\\
\lambda_\infty^\eps-\lambda_{ n}^\eps
\\
\ema \r\Vert_{Y_0^\a(\u_*)}\\
=
{}&\l\Vert \bma
\sum_{i=0}^{\infty} \fr{\deps^i\t_{i}^0}{i!} \eps^i-\sum_{i=0}^{\infty} \fr{\deps^i\t_{n}^0}{i!}\eps^i\\
\sum_{i=0}^{\infty} \fr{\deps^i\p_{i}^0}{i!} \eps^i -\sum_{i=0}^{\infty} \fr{\deps^i\p_{n}^0}{i!}\eps^i\\
\sum_{i=0}^{\infty} \fr{\deps^i\lambda_{i}^0}{i!} \eps^i-\sum_{i=0}^{\infty} \fr{\deps^i\lambda_{n}^0}{i!}\eps^i
\\
\ema \r\Vert_{Y_0^\a(\u_*)} 
\ee
\be
=
{}&\l\Vert \bma
\sum_{i=n}^{\infty} \fr{\deps^i\t_{i}^0}{i!} \eps^i-\sum_{i=n}^{\infty} \fr{\deps^i\t_{n}^0}{i!}\eps^i\\
\sum_{i=n}^{\infty} \fr{\deps^i\p_{i}^0}{i!} \eps^i -\sum_{i=n}^{\infty} \fr{\deps^i\p_{n}^0}{i!}\eps^i\\
\sum_{i=n}^{\infty} \fr{\deps^i\lambda_{i}^0}{i!} \eps^i-\sum_{i=n}^{\infty} \fr{\deps^i\lambda_{n}^0}{i!}\eps^i
\ema \r\Vert_{Y_0^\a(\u_*)} \\
\le {}&
2 \sum_{i=n}^{\infty} \fr{ C^{2i-3}}{i(i-1)}{ \eps}^i \,.
\ee
The claim follows since  
\be
\du(\t_\infty^\eps,\p_\infty^\eps,\lambda_\infty^\eps) = {}&\l(\sum_{i=0}^{\infty} \fr{\du\deps^i\t_{i}^0}{i!}\eps^i,\sum_{i=0}^{\infty} \fr{\du\deps^i\p_{i}^0}{i!}\eps^i,\sum_{i=0}^{\infty} \fr{\du\deps^i\lambda_{i}^0}{i!}\eps^i\r) 
\ee
in $Y_{0}^\a (u_*) $ due to \cref{le limit} and \cref{derivatives estimate}. 
\epr

\section{Proof of Theorem 2.2}\la{se: Main Results Proof}
We apply \cref{thimplicitfunctionIT1 alpha}
to a specific $\ti F$ which is defined below. 

\bde\la{de cutoff}
Let $F,\xi_s$ and $\Xi$ be from \cref{maintheorem}. 
We set
$\ti F(\eps,\xi,x):=\F \chi(\xi),$ 
where $\chi$ is a smooth cutoff function
with $\chi (\xi)=1$ for $|\xi|\le \Xi$ and $\chi (\xi)=0$ for $|\xi|\ge \Xi+1$.
\ede
\noindent
The next lemma follows immediately from the assumptions on $F$ in \cref{maintheorem}.
\ble \la{le: Ftilde vs F}
Let $F$, $\Xi$ be from \cref{maintheorem} and let $ \ti F$ be from \cref{de cutoff}.
Then it holds that
\begin{itemize}
\item[(a)] $\forall ~(\eps,\xi,x) \in (-1,1)\times \l[-\Xi,\Xi \r] \times \R: \ti F(\eps,\xi,x)=\F$\,;
\item[(b)] $\ti F$ satisfies the assumptions  of \cref{thimplicitfunctionIT1 alpha}.  
\end{itemize}
\ele
\noindent
We solve iteratively the equations in 
\cref{thimplicitfunctionIT1 alpha} with the specific $\ti F(\eps,\xi,x):=\F \chi(\xi)$ from \cref{de cutoff} (\cref{thimplicitfunctionIT1 alpha} is applicable due to \cref{le: Ftilde vs F})
and obtain a sequence of solutions, 
which converges due to \cref{le limit}. 
From now on we denote its limit by $(\t_{\infty}^\eps ,\p_{\infty}^\eps ,\lambda_{\infty}^\eps)
$.  
The function $(\t,\p)$ given by \re{form} with
$\bar\xi,\bar u$ satisfying \re{exactODE virtual1},
solves the Cauchy problem \re{SGE1} due to \cref{le limitsolution} and \cref{le: Ftilde vs F}. The claim for $|u_s|\le \tilde C\eps^{}$ follows 
by using \re{exactODE virtual1} and the fundamental theorem of calculus (analogous to the proof of \cite[Lemma 9.2]{Mashkin}). 
\qed
%

\end{document}